\documentclass[10pt,twoside]{article}

\usepackage[utf8]{inputenc}
\usepackage{amsfonts}
\usepackage{amsmath}
\usepackage{amssymb}
\usepackage{amsthm}
\usepackage{amscd}
\usepackage{bbm}
\usepackage{enumerate}

\newcommand{\hocolim}{\operatorname{hocolim}}
\newcommand{\SH}{\mathbf{SH}}
\newcommand{\DM}{\mathbf{DM}}

\newcommand{\QM}{\mathbf{QM}}
\newcommand{\DQM}{\mathbf{DQM}}
\newcommand{\DATM}{\mathbf{DATM}}
\newcommand{\FF}{\mathbb{F}}
\newcommand{\ZZ}{\mathbb{Z}}
\newcommand{\CC}{\mathbb{C}}
\newcommand{\tunit}{\mathbbm{1}}
\newcommand{\Hom}{\operatorname{Hom}}
\newcommand{\id}{\operatorname{id}}
\newcommand{\tr}{{tr}}
\newcommand{\Aone}{{\mathbb{A}^1}}
\newcommand{\eff}{{\text{eff}}}
\newcommand{\colim}{\operatorname{colim}}
\newcommand{\Gm}{\mathbb{G}_m}
\newcommand{\ul}{\underline}

\newcommand{\NN}{\mathbb{N}}
\newcommand{\wequi}{\simeq}
\newcommand{\iRHom}{R\underline{Hom}}

\def\YEAR{\year}\newcount\VOL\VOL=\YEAR\advance\VOL by-1995
\def\firstpage{1}\def\lastpage{1000}
\def\received{}\def\revised{}
\def\communicated{}

\makeatletter
\def\magnification{\afterassignment\m@g\count@}
\def\m@g{\mag=\count@\hsize6.5truein\vsize8.9truein\dimen\footins8truein}
\makeatother

\font\eightrm=cmr8
\font\caps=cmcsc10                    
\font\Caps=cmcsc10 scaled \magstep1   

\makeatletter
\@specialpagefalse\headheight=8.5pt
\def\DocMath{}
\renewcommand{\@evenhead}{%
    \ifnum\thepage>\lastpage\rlap{\thepage}\hfill%
    \else\rlap{\thepage}\slshape\leftmark\hfill{\caps\SAuthor}\hfill\fi}%
\renewcommand{\@oddhead}{%
    \ifnum\thepage=\firstpage{\DocMath\hfill\llap{\thepage}}%
    \else{\slshape\rightmark}\hfill{\caps\STitle}\hfill\llap{\thepage}\fi}%
\makeatother

\def\TSkip{\bigskip}
\newbox\TheTitle{\obeylines\gdef\GetTitle #1
\ShortTitle  #2
\SubTitle    #3
\Author      #4
\ShortAuthor #5
\EndTitle
{\setbox\TheTitle=\vbox{\baselineskip=20pt\let\par=\cr\obeylines%
\halign{\centerline{\Caps##}\cr\noalign{\medskip}\cr#1\cr}}%
	\copy\TheTitle\TSkip\TSkip%
\def\next{#2}\ifx\next\empty\gdef\STitle{#1}\else\gdef\STitle{#2}\fi%
\def\next{#3}\ifx\next\empty%
    \else\setbox\TheTitle=\vbox{\baselineskip=20pt\let\par=\cr\obeylines%
    \halign{\centerline{\caps##} #3\cr}}\copy\TheTitle\TSkip\TSkip\fi%
\centerline{\caps #4}\TSkip\TSkip%
\def\next{#5}\ifx\next\empty\gdef\SAuthor{#4}\else\gdef\SAuthor{#5}\fi%
\ifx\received\empty\relax
    \else\centerline{\eightrm Received: \received}\fi%
\ifx\revised\empty\TSkip%
    \else\centerline{\eightrm Revised: \revised}\TSkip\fi%
\ifx\communicated\empty\relax
    \else\centerline{\eightrm Communicated by \communicated}\fi\TSkip\TSkip%
\catcode'015=5}}\def\Title{\obeylines\GetTitle}
\def\Abstract{\begingroup\narrower
    \parskip=\medskipamount\parindent=0pt{\caps Abstract. }}
\def\EndAbstract{\par\endgroup\TSkip}

\long\def\MSC#1\EndMSC{\def\arg{#1}\ifx\arg\empty\relax\else
     {\par\narrower\noindent%
     2000 Mathematics Subject Classification: #1\par}\fi}

\long\def\KEY#1\EndKEY{\def\arg{#1}\ifx\arg\empty\relax\else
	{\par\narrower\noindent Keywords and Phrases: #1\par}\fi\TSkip}

\newbox\TheAdd\def\Addresses{\vfill\copy\TheAdd\vfill
    \ifodd\number\lastpage\vfill\eject\phantom{.}\vfill\eject\fi}
{\obeylines\gdef\GetAddress #1
\Address #2 
\Address #3
\Address #4
\EndAddress
{\def\xs{4.3truecm}\parindent=0pt
\setbox0=\vtop{{\obeylines\hsize=\xs#1\par}}\def\next{#2}
\ifx\next\empty 
     \setbox\TheAdd=\hbox to\hsize{\hfill\copy0\hfill}
\else\setbox1=\vtop{{\obeylines\hsize=\xs#2\par}}\def\next{#3}
\ifx\next\empty 
     \setbox\TheAdd=\hbox to\hsize{\hfill\copy0\hfill\copy1\hfill}
\else\setbox2=\vtop{{\obeylines\hsize=\xs#3\par}}\def\next{#4}
\ifx\next\empty\ 
     \setbox\TheAdd=\vtop{\hbox to\hsize{\hfill\copy0\hfill\copy1\hfill}
                \vskip20pt\hbox to\hsize{\hfill\copy2\hfill}}
\else\setbox3=\vtop{{\obeylines\hsize=\xs#4\par}}
     \setbox\TheAdd=\vtop{\hbox to\hsize{\hfill\copy0\hfill\copy1\hfill}
	        \vskip20pt\hbox to\hsize{\hfill\copy2\hfill\copy3\hfill}}
\fi\fi\fi\catcode'015=5}}\gdef\Address{\obeylines\GetAddress}

\begin{document}
\Title On the Invertibility of Motives of Affine Quadrics
\ShortTitle 
\SubTitle   
\Author Tom Bachmann
\ShortAuthor 
\EndTitle
\Abstract 
We show that the reduced motive of a smooth affine quadric is invertible as an
object of the triangulated category of motives $\DM(k, \ZZ[1/e])$ (where $k$ is a perfect
field of exponential characteristic $e$). We also establish a
motivic version of the conjectures of Po Hu on products of certain affine
Pfister quadrics. Both of these results are obtained by studying a novel
conservative functor on (a subcategory of) $\DM(k, \ZZ[1/e])$,
the construction of which constitutes the main part of this work.
\EndAbstract
\MSC 
Primary 14C15; Secondary 11E04, 14C25
\EndMSC
\KEY 
Quadric, motive, invertible
\EndKEY
\Address 
Mathematisches Institut
Ludwig-Maximilians-Universität München
Theresienstr. 39
80333 München
Germany
tom.bachmann@zoho.com
\Address
\Address
\Address
\EndAddress

\newtheorem*{theorem*}{Theorem}
\newtheorem*{definition}{Definition}
\newtheorem{theorem}{Theorem}
\newtheorem{proposition}[theorem]{Proposition}
\newtheorem{corollary}[theorem]{Corollary}
\newtheorem{lemma}[theorem]{Lemma}

\numberwithin{equation}{section}


\section{Introduction}
\label{sec:introduction}

For a perfect field $k$, Voevodsky has constructed a triangulated category $\DM(k)$
containing the classical category $Chow(k)$ of Chow motives
\cite{voevodsky-triang-motives}.
Like $Chow(k)$, $\DM(k)$ is a tensor category. We denote the tensor product by
$\otimes = \otimes_{\DM(k)}$ and the unit by $\tunit = \tunit_{\DM(k)}$. As in
any tensor category, we have the notion of \emph{invertible objects}: an object
$E \in \DM(k)$ is called invertible if there exists an object $F \in \DM(k)$ and
an isomorphism $E \otimes F \approx \tunit$. The set of isomorphism classes of
invertible objects forms an abelian group under $\otimes$ and
is called the \emph{Picard group}. We denote it $Pic(\DM(k))$.

Po Hu \cite{HuPicard} was the first to construct interesting elements in
$Pic(\DM(k))$, related to certain low-dimensional quadrics.
In this direction we prove the following result
(see Theorem \ref{thm:invertibility} in Section \ref{sec:app1}).

\begin{theorem*} Let $k$ be a perfect field of exponential characteristic $e$ not two,
$\phi(t_1, \dots, t_n)$ a non-degenerate quadratic form over $k$ and $a \in k^\times$. Write
$X_\phi^a$ for the affine quadric defined by the equation $\phi(t_1, \dots, t_n)
= a$.

Then the reduced motive $\tilde{M}(X_\phi^a) \in \DM(k, \ZZ[1/e])$ is invertible.
\end{theorem*}

This result has a number of predecessors. Work of Voevodsky
can be used to show that
reduced versions of the Rost motives \cite{rost-pfister} are invertible. As
observed by Hu-Kriz \cite[Proposition 5.5]{HuRemarks}, the reduced Rost motives
are reduced motives of affine Pfister quadrics. They go further and
explore analogies with the Hopf invariant one problem. In
\cite{HuPicard} this culminates in certain conjectures about smash products of
affine Pfister quadrics implying their invertibility. Moreover the conjectures
are proven in low dimensions.

The best method the author
knows of attacking the study of Picard groups of tensor categories (to
the extent that it even deserves the name ``method'') is to construct
``realisation functors'' $F: \DM(k) \to \mathcal{C}.$ If $F$ is a tensor
functor, it induces a homomorphism $Pic(\DM(k)) \to Pic(\mathcal{C}).$ If $F$ is
sufficiently nice, and $\mathcal{C}$ sufficiently simple, one may hope to
compute $Pic(\mathcal{C})$ and relate it to $Pic(\DM(k)).$ We mention in passing
that a good test for the ``niceness'' of $F$ seems to be \emph{conservativity}
(i.e. the property that $F$ detects isomorphisms). This will be illustrated
later.

There are well known realisation functors out of $\DM(k)$, but none of them seem
helpful to our problem. If $k \subset \CC$ there is the Hodge realisation, but
this factors through the natural functor $\DM(k) \to \DM(\CC)$ and hence
provides no interesting information about quadrics (which over $\CC$ are
distinguished by only their dimension). There is also étale
realisation, but this factors through $\DM(k) \to \DM_{et}(k)$. In $\DM_{et}(k)$
our problem turns out to be very simple and not indicative of the complexity
encountered in $\DM(k)$ (i.e. in the Nisnevich topology). What we propose in
this work is to construct purpose-built realisation functors $\DM(k) \to
\mathcal{C}$ into big but easy to understand categories. (Actually we do not
quite achieve this; limitations will be explained later.) To motivate our
constructions, we explain two analogous but simpler problems obtained by
replacing $\DM(k)$ by another category.

First let $G$ be a finite group. There exists the \emph{stable $G$-equivariant
homotopy category} $\SH(G)$ \cite{lewis1986equivariant}. Its objects (called genuine $G$-spectra)
are roughly pointed $G$-spaces, where
maps inducing weak equivalences on all fixed point sets have been turned into
isomorphisms, and all representation spheres are invertible objects. If $H \le
G$ is a subgroup, the set of cosets $G/H$ can naturally be turned into a pointed
$G$-space (adding a base point $*$ with trivial action). We denote the
associated spectrum by $\Sigma^\infty G/H_+.$ The objects $\Sigma^\infty G/H_+$
generate $\SH(G).$ There is a functor, called \emph{geometric fixed
points functor}, and denoted $\Phi = \Phi^G: \SH(G) \to \SH$ (where $\SH =
\SH(\{e\})$ is the classical stable homotopy category) which turns out to
be very useful. It is a tensor functor with the property that
$\Phi^G(\Sigma^\infty G/G_+) = S$ (the sphere spectrum), whereas
$\Phi^G(\Sigma^\infty G/H_+) = 0$ for any proper subgroup $H < G$. There are
also natural functors $\SH(G) \to \SH(H)$ (treating $G$-spaces as $H$-spaces)
allowing us to construct the more general geometric fixed points functors
$\Phi^H: \SH(G) \to \SH(H) \to \SH$. As it turns out the collection
$\{\Phi^H\}_H$ (with $H$ ranging over all subgroups of $G$) is as nice as one
may ask (in particular conservative). Consequently these functors were used in
\cite{fausk2001picard} to study $Pic(\SH(G)).$

We now come to a second, more algebraic, example. Let $R$ be a (commutative
unital) ring. Suppose we want to study $Pic(D(R)),$ the Picard group of the
derived category of $R$-modules. Let $m$ be a maximal ideal of $R.$ Recalling
that $D(R)$ can be identified with a subcategory of $K(P(R)),$ the homotopy
category of chain complexes of projective $R$-modules, it is easy to construct a
functor $\Phi^m: D(R) \to D(R/m)$ with the propetry that $\Phi^m(R[0]) =
R/m[0].$ (This is just $\otimes^L_R R/m$.) It turns out that the collection
$\{\Phi^m\}_m$ (where $m$ ranges over all maximal ideals) is as nice as one
needs (at least when restricted to subcategories of sufficiently small objects
in $D(R)$). Moreover the categories $D(R/m)$ are easy to understand.
Consequently, these functors have implicitly been used by Fausk in his study of
the Picard group of derived categories \cite{Fausk-picard-derived}.

Our construction for $\DM^{gm}(k)$ uses a conglomerate of these ideas. The technical
notion of weight structures is the glue that holds our constructions together. We proceed
roughly as follows. Recall that $\DM^{gm}(k)$ is generated as a triangulated category
by the Chow motives. Let $\mathcal{S}$ be the triangulated subcategory generated
by those Chow motives not affording a (non-vanishing) Tate summand. The basic idea is to
consider the (Verdier Quotient) functor $\varphi^k_0: \DM^{gm}(k) \to \DM^{gm}(k) / \mathcal{S}$. The right
hand side does not seem initially easier to understand, but it is at least clear
that it is generated by the images of Tate motives. Using weight structure
theory one obtains a
functor $t: \DM^{gm}(k) / \mathcal{S} \to K^b(Tate)$, where $Tate$ is the
category of (pure) Tate motives, and $K^b$ means bounded chain homotopy category.\footnote{Actually for this to be true we need
to consider $\DM^{gm}(k, \FF)$ where $\FF$ is a \emph{field}. This is not really a
problem.} Combined with base change to arbitrary fields,
we thus obtain a collection of functors $\Phi^l: \DM^{gm}(k) \to \DM^{gm}(l) \to
\DM^{gm}(l)/\mathcal{S} \to K^b(Tate)$.
We not that if $T \in Tate$ is a Tate motive then $\Phi^k(T) = T$. If instead $M
\in Chow$ affords no (non-zero) Tate summands, then $\Phi^k(M) = 0$. This is
rather similar to the geometric fixed points functor $\Phi^G$ from stable
equivariant homotopy theory. Since the general $\Phi^l$ are obtained from
$\Phi^k$ by base change, just as $\Phi^H$ is obtained from the $\Phi^G$
construction by base change (restriction to a subgroup), we will call the
functors $\Phi^l$ ``generalized geometric fixed points functors.''

A natural question is when these functors have
good properties. For our purposes we definitely need \emph{tensor} functors,
which is to say we need $\mathcal{S}$ to be a tensor ideal. This is just not
true in general. However, if instead of looking at the full $\DM^{gm}(k)$ we look at
the subcategory $\DQM^{gm}(k)$ generated by the (products of) smooth projective
quadrics, and use coefficients modulo two, then we can show that $\mathcal{S}$
is a tensor ideal. Moreover, using more properties of weight structures, we
prove the collection of generalized fixed points functors
to be conservative and Pic-injective (i.e. inducing an injection on Picard
groups; see Theorem \ref{thm:final-functors} in Section \ref{sec:app1}):
\begin{theorem*} Let $k$ be a perfect field of exponential characteristic $e$ not two,
and $\Phi^l: \DQM^{gm}(k, \FF_2) \to Tate(\FF_2)$ the functors constructed above.

Then $\{\Phi^l\}_l$, as $l$ ranges over finitely generated extensions of $k$,
forms a conservative, Pic-injective family of tensor triangulated functors.
\end{theorem*}

It is then not hard to use general properties of base change and change of
coefficients for $\DM$ to build a conservative and Pic-injective family for
$\DQM^{gm}(k, \ZZ[1/e])$. It turns out that one additional functor $\Psi: \DQM(k,
\ZZ[1/e]) \to Tate(\ZZ[1/e])$ suffices. (It is related to geometric base change.)

In more detail, the paper is organised as follows. In Section
\ref{sec:chow-motives} we introduce our notations regarding Chow motives and
collect some results. The main idea is to use the absence of degree one
zero-cycles in a variety to conclude that it is free of Tate summands in a
strong sense. This observation is what will allow us in a later section to
establish that our ``geometric fixed points functors'' $\Phi^l$ are tensor.

In Section \ref{sec:review-vvm} we review in some detail the categories $\DM(k,
A)$ (triangulated motives over the perfect field $k$ with coefficients in the
commutative ring $A$) and their behaviour
under change of coefficients and base. All the material is well known, but
sometimes hard to source. We then construct a convenient conservative and Pic-injective
collection on $\DM(k, A)$. The targets are always $\DM(k', A')$ with either
$k$ simplified (e.g. $k'$ separably closed) or $A$ simplified (e.g. $A'$ a
field).

Section \ref{sec:weights} constitutes the technical heart of our work. We
first rapidly review Bondarko's theory of weight structures. After that we carry
out the programme outlined above, of constructing a conservative and
Pic-injective family of functors $\{\Phi^l\}_l: \DQM^{gm}(k, \FF_2) \to
K^b(Tate(\FF_2))$.

The remaining sections contain applications.
In Section \ref{sec:app1} we prove that
all affine quadrics have invertible motives. This is rather satisfying,
since affine quadrics are fairly natural ``generalised spheres.'' Also the
result has been known in the étale topology for a long time. Compare the
beginning of this introduction for a history of this problem.

Section \ref{sec:app2} contains the second set of applications. In
\cite[Conjecture 1.4]{HuPicard} Po Hu has stated certain conjectures about the
motivic spectra of affine Pfister quadrics, namely certain formulas they should
satisfy under smash product. We establish the analogues (or ``images'') of these
formulas in $\DM(k)$ by an easy computation involving our fixed points functors.

The list of applications of our methods does not end here, but the amount of
material we want to stuff into one article does. As directions of future work,
let us mention the following possibilities. The structure of $Pic(\DQM(k))$ can
be investigated. One may replace
the set of projective quadrics by projective homogeneous varieties for a fixed
group $G.$ Also using (almost) the same methods, it is possible to study $\DATM(k),$
the subcategory of $\DM(k)$ generated by $M(Spec(l))\{i\}$ for $l/k$ finite
separable and $i \in \ZZ,$ i.e. Artin-Tate motives. This will be treated in
forthcoming work.

We also note that our results for $\DM$ have applications to the study of the
stable motivic homotopy category $\SH(k)$. In forthcoming work
\cite{bachmann-hurewicz} we
show that if $k$ is a field of finite 2-étale cohomological dimension, then the
functor $\SH(k) \to \DM(k)$ is conservative and Pic-injective, when restricted
to compact spectra. Consequently the suspension spectral of affine quadrics are
also invertible, and the Hu conjectures hold for spectra (over such fields).

Whenever we talk about quadrics or quadratic forms, we shall assume that the base field has
characteristic different from two. This will be restated with the most important
theorems.

Our results are stated over perfect base fields, because this is when $\DM(k)$
is best understood. However actually everything goes through over arbitrary base
fields, using \cite{integral-mixed-motives}.

Throughout this text, we will omit brackets around the arguments of functors
whenever convenient. For example $MX$ means the same thing as $M(X)$.

The author wishes to thank Fabien Morel for suggesting this topic of
investigation and for providing many helpful insights,
and Mikhail Bondarko for comments on a draft of this paper. He also wishes to
thank an anonymous referee for many helpful comments and suggestions.

\section{Some Results about Chow Motives}
\label{sec:chow-motives}

We begin with some notation. We take for granted the notion of an \emph{additive
category}. An additive category $\mathcal{C}$ is called \emph{Karoubi-closed} if
every idempotent endomorphism of an object of $\mathcal{C}$ corresponds to a
direct sum decomposition. By a \emph{tensor category} we mean an additive
category provided with a suitably compatible symmetric monoidal
structure \cite[Section 1]{tannakian-cats}. In particular this means that the monoidal operation is
bi-additive. We shall always denote the monoidal operation by $\otimes =
\otimes_\mathcal{C}$ and call it tensor product. The tensor unit is generically
denoted $\tunit = \tunit_\mathcal{C}.$

Now our conventions regarding Chow motives. By $SmProj(k)$ we denote the
category of smooth projective varieties over the field $k$. It is a symmetric
monoidal category using cartesian product as monoidal product. We shall assume
understood the existence and functoriality properties of the \emph{Chow ring}
$A^*(X).$ Grading is by codimension and the equivalence relation we use is
rational equivalence. Lower index means grading by dimension. For convenience if
$\FF$ is any coefficient ring, we put $A^*(X, \FF) = A^*(X) \otimes_\ZZ \FF.$ It
is then possible to construct a Karoubi-closed tensor category $Chow(k, \FF)$
together with a covariant symmetric monoidal functor $M = M_\FF: SmProj(k) \to
Chow(k, \FF)$ which has the following properties. The unit object is
$\tunit_{Chow(k, \FF)} = \tunit = M(Spec(k))$. There exists an object
$\tunit\{1\}$ such that $M(\mathbb{P}^1) \approx \tunit \oplus \tunit\{1\}.$ We
call $\tunit\{1\}$ the Lefschetz motive. It is invertible. For any $n \in \ZZ$
and $M \in Chow(k, \FF)$ we write $M\{n\} := M \otimes \tunit\{1\}^{\otimes n}$.
For any $X, Y \in SmProj(k)$ and $i, j \in \ZZ$ we have
\[ \Hom_{Chow(k, \FF)}(M(X)\{i\}, M(Y)\{j\}) = A_{\dim{X} + i - j}(X \times Y).  \]
 In particular we have
$\Hom(MX, \tunit\{i\}) = A^i(X, \FF)$ and $\Hom(\tunit\{i\}, MX) = A_i(X, \FF).$
Composition is by the usual push-pull convolution.

In the remainder of this section we collect some results about Chow motives
which we will need throughout the article. None of them are hard, so probably
most of this is well known.

Recall first that if $l/k$ is a field extension then the functor $SmProj(k) \to SmProj(l), X
\mapsto X_l$ induces a functor $Chow(k, \FF) \to Chow(l, \FF)$ called \emph{base
change} and denoted $M \mapsto M_l$. We need to know something about this in the
inseparable case.

\begin{lemma} \label{lemm:chow-insep-basechange}
Let $l/k$ be a purely inseparable extension of fields of characteristic $p$ and
$\FF$ a coefficient ring in which (the image of) $p$ is invertible. Then the base change
$Chow(k, \FF) \to Chow(l, \FF)$ is fully faithful.
\end{lemma}
\begin{proof}
It suffices to prove that for $X \in SmProj(k)$ we have $A_*(X, \FF) =
A_*(X_l, \FF)$. By the definition of rational equivalence as in \cite[Section
1.6]{fulton-intersection} it is enough to show that $Z_*(X, \FF) \to Z_*(X_l,
\FF)$ is an isomorphism for all $X$.

Let $Z \subset X$ be a reduced closed subscheme and $|Z_l|$ the reduced closed
subscheme underlying $Z_l$. Then the image of $[Z]$ under $Z_*(X, \FF) \to
Z_*(X_l, \FF)$ is $n [|Z_l|]$, where $n$ is the multiplicity of $Z_l$. This is
easily seen to be a power of $p$, whence $Z_*(X, \FF) \to Z_*(X_l, \FF)$ is
injective. It is also surjective since $X_l \to X$ is a homeomerphism on
underlying spaces. This concludes the proof.
\end{proof}

We now investigate ``Tate summands''. Denote by $Tate(k, \FF) \subset Chow(k,
\FF)$ the smallest (strictly) full Karoubi-closed additive subcategory containing
$\tunit\{i\}$ for all $i$. This is independent up to equivalence of $k$ and we
will just write $Tate(\FF)$ if no confusion can arise. (It is a tensor category.)

We say $M \in Chow(k, \FF)$ is \emph{Tate-free} if whenever $M \approx T \oplus
M'$ with $T \in Tate(k, \FF)$, then $T \approx 0$. The next proposition holds in
much greater generality, but this version is all we need.

\begin{proposition} \label{prop:tatefree-decomp}
Let $\FF$ be a finite ring and $M \in Chow(k, \FF)$. Then there exist $T \in
Tate(\FF)$ and $M' \in Chow(k, \FF)$ with $M'$ Tate-free and $M \approx T \oplus
M'$.
\end{proposition}
\begin{proof}Splitting off
Tate summands inductively, the only problem which could occur is that $M$
might afford arbitrarily large Tate summands. The impossibility of this follows
(for example) from the finiteness of étale cohomology of complete varieties
\cite[Corollary VI.2.8]{MilneEtaleCohomology}.
\end{proof}

\begin{lemma} \label{lemm:sums-tatefree}
Let $\FF$ be a field. Then if $M, N \in Chow(k, \FF)$ are Tate-free so is $M
\oplus N$.
\end{lemma}
\begin{proof}
A motive with $\FF$-coefficients is Tate-free if and only if it affords no
summand of the form $\tunit\{n\}$ for any $n$.

Let $i: \tunit\{n\} \to M \oplus
N$ and $p: M \oplus N \to \tunit\{n\}$ be inclusion of and projection to a
summand, for $M, N$ arbitrary. Write $i = (i_M, i_N)^T$ and $p = (p_M, p_N)$. Then $\id = pi = p_M i_M
+ p_N i_N$. Since $\Hom(\tunit\{n\}, \tunit\{n\}) = \FF \ne 0$ we must have $p_M i_M
\ne 0$ or $p_N i_N \ne 0$. Suppose the former holds. Then since $\FF$ is a field
we may replace $i_M$ by a multiple $c i_M$ such that $p_M (c i_M) = 1$. Thus
$\tunit\{n\}$ is a summand of $M$. Similarly in the other case. This establishes
the contrapositive of the lemma.
\end{proof}

\begin{lemma} \label{lemm:tatefree-factoring}
Let $\FF$ be a field. Then any morphism in $Tate(k, \FF)$ factoring through a
Tate-free object is zero.
\end{lemma}
\begin{proof}
Since $\FF$ is a field any Tate motive is a sum of $\tunit\{n\}$ for various
$n$, so it suffices to consider a morphism $\tunit\{n\} \to \tunit\{m\}$
factoring through a Tate-free object. Since $\Hom(\tunit\{n\}, \tunit\{m\}) = 0$
for $n \ne m$ we may assume $n = m$. Consider $a \in \Hom(\tunit\{n\}, M)$ and
$b \in \Hom(M, \tunit\{n\})$. If $ba \ne 0$ then there exists $c \in \FF$
such that $(cb) a = \id$. It follows that $(cb), a$ present $\tunit\{n\}$ as a
summand of $M$. This establishes the contrapositive.
\end{proof}

We need tools to recognise Tate-free motives. To do so, we introduce some
more notation. For $X \in SmProj(k)$ there exists the degree map $deg: A_0(X,
\FF) \to \FF$ (corresponding to pushforward along the structure map $\Hom(\tunit,
MX) \to \Hom(\tunit, \tunit)$). Write $I_\FF(X) = deg(A_0(X, \FF))$ for the
image of the degree map. This is the ideal inside $\FF$ generated by the degrees
of closed points. The utility of this
notion is as follows.

\begin{lemma} \label{lemm:IF-tatefree}
Let $\FF$ be a field and suppose $I_\FF(X) \ne \FF$. Then $MX$ is Tate-free.
\end{lemma}
\begin{proof}
As before $MX$ is Tate-free if and only if it affords no summand $\tunit\{N\}$
for any $N$.
Given $i \in \Hom(\tunit\{N\}, MX) = A_N(X, \FF)$ and $p \in \Hom(MX,
\tunit\{N\}) = A^N(X, \FF)$, the composite $pi \in \Hom(\tunit\{N\},
\tunit\{N\}) = \FF$ is obtained by push-pull convolution. In this case it is
just $deg(p \cap i)$ and so is contained in $I_\FF(X).$ Thus $pi \ne 1$ and $(p,
i)$ is \emph{not} a presentation of $\tunit\{N\}$ as a summand of $X$.
\end{proof}

\begin{lemma} \label{lemm:IF-products}
Let $X, Y \in SmProj(k)$. Then $I_\FF(X \times Y) \subset I_\FF(X) \cap
I_\FF(Y)$.
\end{lemma}
\begin{proof}
We recall that $I_\FF(X \times Y)$ is just the ideal generated by degrees of
closed points. So let $z \in X \times Y$ be a closed point. Then $z \to X \times
Y$ corresponds to morphisms $z \to X$ and $z \to Y$. It follows that $deg(z) \in
I_\FF(X)$ and similarly $deg(z) \in I_\FF(Y)$. This implies the result.
\end{proof}

Suppose $S \subset SmProj(k)$ is a set of smooth projective varieties. We write
$\langle S \rangle^{\otimes, T}_{Chow(k, \FF)}$ for the smallest strictly full, additive,
Karoubi-closed, tensor subcategory of $Chow(k, \FF)$ containing all Tate motives
and also $MX$ for each $X \in S$. Assuming $\FF$ is a field, this means that a
general object of $\langle S \rangle^{\otimes, T}_{Chow(k, \FF)}$ is (isomorphic
to) a summand of
\[ T \oplus M(X_1^{(1)} \times \dots \times X_{n_1}^{(1)})\{i_1\} \oplus \dots
     \oplus M(X_1^{(m)} \times \dots \times X_{n_m}^{(m)})\{i_m\}, \]
with $T \in Tate(\FF)$ and $X_i^{(j)} \in S, i_r \in \ZZ$.

The following proposition (or rather its failure to generalise)
is the basic reason why in the construction of fixed point
functors we will need to restrict to subcategories.

\begin{proposition} \label{prop:tatefree-presentation}
Let $\FF$ be a finite field and $S \subset SmProj(k)$ be such that $I_\FF(X) = 0$ for
all $X \in S$ (i.e. such that all closed points of $X$ have degree divisible by
the characteristic of $\FF$). Then any object $M \in \langle S \rangle^{\otimes,
T}_{Chow(k, \FF)}$ can be written as $T \oplus M'$, where $T \in Tate(\FF)$ and
$M'$ is (isomorphic to) a summand of 
\[ M(X_1^{(1)} \times \dots \times X_{n_1}^{(1)})\{i_1\} \oplus \dots
     \oplus M(X_1^{(m)} \times \dots \times X_{n_m}^{(m)})\{i_m\}, \]
for some $X_i^{(j)} \in S, i_r \in \ZZ$. Moreover any such $M'$ is Tate-free.
\end{proposition}
\begin{proof}
By Lemma \ref{lemm:IF-products} we know that $I_\FF(X_1^{(j)} \times \dots
X_{n_j}^{(j)}) = 0$ and thus by Lemmas \ref{lemm:IF-tatefree} and
\ref{lemm:sums-tatefree} we conclude that any $M'$ as displayed is indeed
Tate-free. So it suffices to establish the first part.

By definition we may write
\[ M \oplus M'' \approx T \oplus M(X_1^{(1)} \otimes \dots \otimes X_{n_1}^{(1)})\{i_1\} \oplus \dots
     \oplus M(X_1^{(m)} \otimes \dots \otimes X_{n_m}^{(m)})\{i_m\}, \]
with $T \in Tate(\FF)$ and $X_i^{(j)} \in S$. Using Proposition
\ref{prop:tatefree-decomp} we write $M \oplus M'' \approx M' \oplus M''' \oplus
T'$, where $M', M''$ are maximal Tate-free summands in $M, M''$ respectively and
$T'$ is Tate. Writing out the inverse isomorphisms $M' \oplus M''' \oplus T'
\leftrightarrows T \oplus M(X_1^{(1)} \dots) \oplus \dots$ in matrix form and
using Lemma \ref{lemm:tatefree-factoring} we conclude that $T' \approx T$ via
the induced map. The Lemma below yields that $M' \oplus M'' \approx M(X_1^{(1)}
\dots) \oplus \dots$. This finishes the proof.
\end{proof}

\begin{lemma} \label{lemm:cancel-summand}
Let $\mathcal{C}$ be an additive category and let $U, T, X, T' \in \mathcal{C}$
be four objects. Suppose we are given an isomorphism $\phi: U \oplus T \to X
\oplus T'$ such that the component $T \to T'$ is also an isomorphism. Then there
is an isomorphism $\tilde \phi: U \to X$.
\end{lemma}
\begin{proof}
Let us write
\[ \phi = \begin{pmatrix} \alpha & a \\ b & f \end{pmatrix} \quad
   \psi = \begin{pmatrix} \beta & a' \\ b' & g \end{pmatrix},\]
where $\psi$ is the inverse of $\phi$. By assumption $f$ is an
isomorphism. Writing out $\phi \psi = \id_{X\oplus T}$ and $\psi \phi =
\id_{U\oplus T'}$ one obtains
\[ b \beta = - fb' \quad \beta a = - a'f \quad \alpha \beta + ab' = \id_U \quad
\beta \alpha + a'b = \id_X. \]
Put $\tilde \alpha = \alpha - a f^{-1} b: U \to X$. Then the above relations imply that
$\tilde \alpha$ is an isomorphism with inverse $\beta$.
\end{proof}

\section{Review of Voevodsky Motives}
\label{sec:review-vvm}

In this section we collect some facts about $\DM(k,A)$ that we will need for the
applications. Most of this is well-known, but in some cases we were unable to
locate adequate references.
We will assume throughout that $k$ is a perfect field, and re-state this assumption with
each theorem.

Fix a ring $A$. We follow the construction of $\DM(k,A)$ and $\DM^\eff(k,A)$ in
\cite{cisinski2009local}.

Write $Sm(k)$ for the symmetric monoidal category of smooth schemes over $k$
(monoidal operation being cartesian product) and $Cor(k)$ for the symmetric
monoidal category with same objects as $Sm(k)$ but morphisms given by finite
correspondences. There is a natural monoidal functor $Sm(k) \to
Cor(k)$. Write $Shv^{tr}(k, A)$ for the abelian category of Nisnevich sheaves of
$A$-modules on $Cor(k)$, i.e. (additive) presheaves $Cor(k)^{op} \to A\text{-Mod}$ such that the
restriction $Sm(k)^{op} \to Cor(k) \to A\text{-Mod}$ is a sheaf in the Nisnevich
topology. There is a functor $A_{tr} \bullet: Sm(k) \to Shv^{tr}(k, A)$ sending $X
\in Sm(k)$ to the presheaf with transfers it represents (which turns out to be a
sheaf). The category $Shv^{tr}(k, A)$ carries a right exact tensor structure making
$A_\tr$ a monoidal functor.

The category $\DM^\eff(k,A)$ is then the $\Aone$-local derived category of
$Shv^{tr}(k,A)$ \cite[Example 3.15]{cisinski2009local}. We write $L_\Aone:
D(Shv^{tr}(k,A)) \to \DM^\eff(k,A) \subset D(Shv^{tr}(k,A))$
for the localisation functor and denote the composite $Sm(k)
\xrightarrow{A_{tr}} Shv^{tr}(k,A) \to \DM^\eff(k,A)$ by $M^\eff_A$ or $M^\eff$
if no confusion can arise. The category $\DM^\eff(k,A)$ is
compactly generated, and the subcategory of compact objects $\DM^{\eff,gm}(k,A)$
is the thick subcategory generated by $M^\eff_A(X)$ for $X \in Sm(k)$
\cite[Example 5.5]{cisinski2009local}. The category $\DM^\eff(k,A)$ carries a
symmetric monoidal structure making $M^\eff_A: Sm(k) \to \DM^\eff(k,A)$ a
monoidal functor \cite[Example 2.4]{cisinski2009local}.

If $X \in Sm(k)$ then we write $\tilde{M}^\eff X$ for the homotopy fibre of
$M^\eff X \to M^\eff Spec(k) = \tunit$. This is the \emph{reduced} (effective)
motive of $X$. If $(X, x)$ is a pointed scheme then
there is a canonical isomorphism $M^\eff X \wequi \tilde{M}^\eff X \oplus
M^\eff x = \tilde{M}^\eff X \oplus \tunit$. In this situation we will write
$M^\eff(X,x)$ for $\tilde{M}^\eff X$, with this direct sum decomposition
understood.

Throughout this text, we write $\Gm$ for the pointed scheme $(\Aone \setminus 0,
1)$.

The next step is to stabilise $\DM^\eff(k,A)$ by inverting $M^\eff \Gm$ in the
monoidal structure. Here we depart slightly from the notation of
\cite{cisinski2009local} and write $Sp(Shv^{tr}(k,A))$ for the abelian category
of symmetric spectra in $Shv^{tr}(k,A)$ \cite[Section 6]{cisinski2009local}.
This is in keeping with the notation in \cite[Section
5.3]{triangulated-mixed-motives}. Then $\DM(k,A)$ is the $\Omega-\Aone$-local
derived category of $Sp(Shv^{tr}(k,A))$ \cite[Example 6.25]{cisinski2009local}.
There is an adjunction $\Sigma^\infty:
Shv^{tr}(k,A) \leftrightarrows Sp(Shv^{tr}(k,A))$ extending to an adjunction
\[ \Sigma^\infty: \DM^\eff(k,A) \leftrightarrows \DM(k,A): \Omega^\infty. \]
We write $M = M_A: Sm(k) \to \DM(k,A)$ for the evident composite, and call $MX$
for $X \in Sm(k)$ the motive of $X$. Similarly for the reduced motive
$\tilde{M}X$. As before,
$\DM(k,A)$ is a compactly generated tensor triangulated category. The
subcategory $\DM^{gm}(k,A)$ of compact objects is the thick triangulated
subcategory generated by the objects of the form $MX \otimes (M\Gm)^{\otimes i}$
for $X \in Sm(k)$ and $i \in \ZZ$.

So far all of this is very formal and the base $k$ did not really enter. Since
$k$ is perfect, we have Voevodsky's remarkable results at hand. Firstly, a
bounded above complex in $D(Shv^{tr}(k,A))$ is $\Aone$-local if and only if its
homology sheaves are \cite[Proposition 14.8]{lecture-notes-mot-cohom}. Moreover
a model for $L_\Aone$ is given by the $\Aone$-chain complex $C_*$
\cite[Corollary 14.9]{lecture-notes-mot-cohom}.

Next there is the
cancellation theorem: for $E, F \in \DM^\eff(k,A)$ we have $\Hom(E \otimes M\Gm,
F \otimes M\Gm) = \Hom(E, F)$ \cite[Corollary 4.10]{voevodsky2002cancellation}.
(Voevodsky only states this for $E, F$ bounded above, but the general case
follows using compact generation and taking limits: Let $\mathcal{C}_1 \subset
\DM^\eff(k,A)$ be the class of objects $F$ such that $\Hom(MX[i], F) = \Hom(MX[i]
\otimes M\Gm, F \otimes M\Gm)$ for all $X \in Sm(k)$. Then $\mathcal{C}_1$ is
closed under cones (by the five lemma), shifts, isomorphisms and arbitrary sums
(because the $MX$ are compact) and
contains all bounded above complexes, hence $\mathcal{C}_1 = \DM^\eff(k,A)$.
Next let $\mathcal{C}_2 \subset \DM^\eff(k,A)$ be the class of objects $E$ such
that $\Hom(E, F) = \Hom(E \otimes M\Gm, F \otimes M\Gm)$ for all $F \in
\DM^\eff(k,A)$. Then $\mathcal{C}_2$ is closed under shifts, cones (by the
five lemma) isomorphisms and arbitrary sums, and contains $\DM^{gm}(k,A)$ (since
$\mathcal{C}_1 = \DM^\eff(k,A)$) and hence $\mathcal{C}_2 = \DM^\eff(k,A)$.)
This implies that
$\Sigma^\infty: \DM^\eff(k,A) \to \DM(k,A)$ is fully faithful. Indeed if $E \in
D(Shv^{tr}(k,A))$ is $\Aone$-local then $\Sigma^\infty E \in
D(Sp(Shv^{tr}(k,A)))$ is an $\Omega$-spectrum by cancellation (i.e. $E \wequi
\iRHom(M\Gm, E \otimes M\Gm)$).

Finally there are the homotopy $t$-structures: the category $\DM^\eff(k,A)$
has a non-degenerate $t$-structure with heart the category of homotopy
invariant sheaves with transfers, and $\DM(k,A)$ affords a non-degenerate
$t$-structure with heart the category of homotopy modules with transfers. The
functor $\Omega^\infty$ is $t$-exact and hence $\Sigma^\infty$ is
right-$t$-exact.

We will now discuss functoriality of $\DM(k,A)$ in $k$ and $A$. The tool to do
this is \cite[Proposition 3.11]{cisinski2009local}, which implies that
if $f^*: Shv^{tr}(k,A) \to Shv^{tr}(l,B)$ is a functor which preserves colimits,
coverings, and multiplication by $\Aone$ and $\Gm$,
then there are induced adjunctions
\[ Lf^*: \DM^\eff(k,A) \leftrightarrows \DM^\eff(l,B): Rf_* \]
\[ Lf^*: \DM(k,A) \leftrightarrows \DM(l,B): Rf_*. \]
Here $Lf^*$ commutes with $\Sigma^\infty$ and $Rf_*$ commutes with
$\Omega^\infty$. If the functor $f^*$ we started with was monoidal then so is the
extended functor $Lf^*$.

We now specialise to base change. For this, let $f: Spec(l) \to Spec(k)$ be an
extension of (perfect) fields. Then $f^*: Shv^{tr}(k,A) \to Shv^{tr}(l,A)$
satisfies the requirements outlined above and so we get $Lf^*: \DM(k,A)
\leftrightarrows \DM(l,A): Rf_*$, and similarly for $\DM^\eff$. If $f$ is finite
separable then $f_\#: Sm(l) \to Sm(k)$ induces $f_\#: Shv^{tr}(l,A) \to
Shv^{tr}(k,A)$ and then $Lf_\#: \DM(l,A) \leftrightarrows \DM(k,A): Rf^*$.
See also \cite[Example 6.25]{cisinski2009local} again. In this situation we have
$Rf^* = f^* = Lf^*$.

The following result is surely well-known, but we could not find a reference, so
include the easy proof.

\begin{proposition} \label{prop:f*-conservative}
Let $f: Spec(l) \to Spec(k)$ be an algebraic (separable) extension of the perfect field
$k$. Then $Lf^*: \DM(k,A) \to \DM(l,A)$ is $t$-exact.

Suppose that $A$ is a ring such that for each finite subextension $l/l'/k$, the (image
of the) integer $[l':k]$ is a unit in $A$.

Then $Lf^*: \DM(k,A) \to \DM(l,A)$ is conservative and $t$-exact.
\end{proposition}
Thus we shall write $Lf^* = f^*$ also in this situation.
\begin{proof}
Since $f^*: Shv^{tr}(k,A) \to Shv^{tr}(l,A)$ preserves homotopy invariant
sheaves, it follows that $Lf^*: D(Sp(Shv^{tr}(k,A))) \to D(Sp(Shv^{tr}(k,A)))$
preserves $\Aone$-local objects. Since $f^*: Shv^{tr}(k,A) \to Shv^{tr}(l,A)$
preserves contractions, $Lf^*: D(Sp(Shv^{tr}(k,A))) \to D(Sp(Shv^{tr}(l,A)))$
preserves $\Aone-\Omega$-local objects. Thus $t$-exactness of 
$Lf^*: D(Sp(Shv^{tr}(k,A))) \to D(Sp(Shv^{tr}(l,A)))$ implies $t$-exactness of
$Lf^*: \DM(k,A) \to \DM(l,A)$.

It thus remains to show: if $F \in Shv^{tr}(k,A)$
and $f^*F = 0$, then $F=0$. Let $X \in Sm(k)$, $x \in F(X)$. It suffices to show
that $x = 0$. By \cite[Proposition II.2.2 and
Lemma II.3.3]{MilneEtaleCohomology} we have $0 = f^*x \in F(X \otimes_k
l) = \colim_{l/l'/k} F(X \otimes_k l')$, where the colimit is over finite
subextensions. Thus there exists a finite subextension
$l/l'/k$ with $(l'/k)^*(x) = 0$. But then by a transfer argument one finds that
$[l':k]x = 0$, whence $x = 0$ since $[l':k]$ is a unit in $A$ by assumption.
\end{proof}

Next we consider change of coefficients. The construction and basic properties
must be well known, but again we could not find convenient references. Let
$\alpha: A \to B$ be a ring homomorphism. There is a natural adjoint functor
pair \[\alpha_\#: Shv^{tr}(k, A) \rightleftarrows Shv^{tr}(k, B): \alpha^*.\]
Here $\alpha_\# F$ is the sheaf associated to $X \mapsto F(X) \otimes_A B$
and $\alpha^*F(X) = F(X)$, viewed as an
$A$-module. In particular $\alpha_\#$ is monoidal and preserves colimits, so
\cite[Proposition 3.11]{cisinski2009local} applies to give us adjunctions
\[ L\alpha_\#: \DM^\eff(k,A) \to \DM^\eff(k,B): R\alpha_\# \]
\[ L\alpha_\#: \DM(k,A) \to \DM(k,B): R\alpha_\#. \]

In order to manipulate these functors efficiently, we need a standard result.

\begin{lemma} \label{lemm:hocolim-decomp}
Any object $E \in \DM(k,A)$ is a filtered homotopy colimit of objects of the
form $E'\otimes \Gm^{\otimes i}$ with $E' \in \DM^\eff(k,A)$ bounded above and $i
\in \ZZ$.
\end{lemma}
\begin{proof}
Let $\mathcal{DM}^\eff(k,A)$ be a model for $\DM(k,A)$. Then we may
alternatively model
$\DM(k,A)$ via $Spt^\Sigma(\mathcal{DM}^\eff(k,A), M\Gm)$, i.e. via (symmetric)
$M\Gm$-spectra in $\mathcal{DM}^\eff(k,A)$.

Let $E \in \DM(k,A) = Ho(Spt^\Sigma(\mathcal{DM}^\eff(k,A), M\Gm))$ have fibrant
replacement $(E_1, E_2, \dots)$. Then $E \wequi \hocolim_i
(E_i)_{\ge -i} \otimes \Gm^{\otimes -i}$, as one sees immediately by computing
the homotopy sheaves on both sides.
\end{proof}
In order to use this, recall that any left adjoint functor of triangulated
categories (assumed to have all countable coproducts) commutes with
filtered homotopy colimits (by filtered we always mean $\omega$-filtered, i.e.
with a countable indexing set) and if it additionally preserves a compact
generating set, then its right adjoint also commutes with filtered homotopy
colimits. Thus essentially all our functors commute with filtered homotopy
colimits. In particular $R\alpha^*$ commutes with filtered homotopy colimits.

\begin{lemma} \label{lemm:Ralpha*-comm}
Let $k$ be perfect and $\alpha: A \to B$ a ring homomorphism. Then the functors
$R\alpha^*$ commute with $\Sigma^\infty$. In fact for $E \in \DM^\eff(k,B)$ we
have $R\alpha^*(E \otimes \Gm) \wequi R\alpha^*(E) \otimes \Gm$.
\end{lemma}
Of course $L\alpha_\#$ always commutes with $\Sigma^\infty$, for formal reasons.
\begin{proof}
Let $E \in \DM^\eff(k,B)$. Then $\Sigma^\infty E \wequi (E, E \otimes \Gm, E \otimes
\Gm^{\otimes 2}, \dots)$. (Here by $E \otimes \Gm$ we mean the derived tensor
product in $\DM^\eff$, in particular this notation implies an $\Aone$-local
object.) By the cancellation theorem, this is an
$\Omega$-spectrum. It follows that $R\alpha^* \Sigma^\infty E = (R\alpha^* E,
R\alpha^* (E \otimes \Gm), \dots)$. It is thus enough to show that $R\alpha^*(E
\otimes \Gm) \wequi R\alpha^*(E) \otimes \Gm$.

Since $L\alpha_\#$ is symmetric monoidal $R\alpha^*$ is lax symmetric monoidal
and there is a natural comparison map. Since $\DM^\eff(k,B)$ is generated as a
localising subcategory by $MX$ for $X \in Sm(k)$ and $\otimes, R\alpha^*$ commute
with arbitrary sums, we may assume $E = MX$. In this case a fibrant model of $E
\otimes \Gm$ is
given by $C_* B_{tr} (X_+ \wedge \Gm)$. Resolving $B$ freely as an $A$-module,
it follows that $R\alpha^*E \in
\DM^\eff(k,A) \otimes \Gm$. A calculation using adjunction and the cancellation
theorem allows us to conclude by the Yoneda lemma.
\end{proof}

\begin{proposition} \label{prop:alpha-flat}
Let $k$ be perfect, $\alpha: A \to B$ be flat, $E \in \DM^{gm}(k,A)$ and $F \in
\DM(k,A)$.
Then \[\Hom(E,F) \otimes_A B \approx \Hom(L\alpha_\# E, L\alpha_\# F). \]
\end{proposition}
\begin{proof}
By Lemma \ref{lemm:hocolim-decomp} (and twisting $E$), we may assume
that $F \in \DM^\eff(k,A)$ and is bounded above. Then by the cancellation theorem we may
assume that $E \in \DM^{gm,\eff}(k,A)$ as well.

There is
a natural map from the left hand side to the right hand side. Using the $5$-lemma and
the fact that $\DM^{gm,\eff}(k,A)$ is generated by $MX$ for
$X \in Sm(k)$, we
may reduce to $E = MX$ (shifting $F$ if necessary). In this case $\Hom(MX, F)$ is given by the
hypercohomology $H^0(X, F^\bullet)$. Since $\otimes_A B$ is exact it commutes
with hypercohomology and preserves $\Aone$-invariance of cohomology sheaves,
so we have $H^0(X, F^\bullet)
\otimes_A B = H^0(X, F^\bullet \otimes_A B) = H^0(X, (L\alpha_\#
 F)^\bullet)$.
\end{proof}

\begin{proposition} \label{prop:bockstein-triangle}
Let $k$ be perfect, $A$ a ring, $a \in A$ a non zero divisor and $\alpha: A \to
A/(a)$ the natural map. Then for $E \in \DM(k,A)$ there is a natural
distinguished triangle
\[ E \xrightarrow{\cdot a} E \to R\alpha^* L\alpha_\# E. \]
\end{proposition}
This triangle yields the typical \emph{Bockstein sequences} one expects for
reduction of coefficients.
\begin{proof}
By Lemmas \ref{lemm:hocolim-decomp} and \ref{lemm:Ralpha*-comm} we may assume
that $E \in \DM^\eff(k,A)$ and is bounded above.

In this case $R\alpha^* L\alpha_\# E$ is computed by resolving
$E$ by a complex of representable sheaves $C^\bullet$ and then $C^\bullet/(a)$
is a model for $R\alpha^* L\alpha_\# E$. (Note that since $C^\bullet$ has
homotopy invariant cohomology, so does $\alpha_\# C^\bullet = C^\bullet/(a)$, by
considering the (ordinary) Bockstein sequence. Hence we may apply $\alpha^*$ immediately to
$\alpha_\# C^\bullet$ instead of having to $\Aone$-localise first.) Since $a$ is
not a zero divisor the sequence $0 \to C^\bullet \to C^\bullet \to C^\bullet/(a)
\to 0$ is exact and yields the desired triangle.
\end{proof}

With this preparation out of the way, we can prove our conservativity and
Pic-injectivity theorem. Recall that $\Hom_{\DM(k,A)}(\tunit, \tunit[i]) = A$ if
$i = 0$ and $= 0$ else.

\begin{theorem} \label{thm:conservative-on-DM}
Let $k$ be a perfect field and $A$ a PID of characteristic zero. Let $f:
Spec(k^s) \to Spec(k)$ be a separable closure.

The collection of functors $\{f^* \} \cup \{L\alpha_{\pi \#}\}_\pi$ is conservative. If
$A$ has primes of arbitrary large characteristic, the collection is also
Pic-injective (both on $\DM(k,A)$). Here $\alpha_\pi: A \to A/(\pi)$ runs through
the primes of $A$.
\end{theorem}
We could prove essentially the same theorem with $A$ replaced by a Dedekind domain (of
characteristic zero) with only slightly more work.
\begin{proof}
We first show conservativity. Let
$E \in \DM(k,A)$ with $L\alpha_{\pi\#}E = 0$ for all $\pi$ and $f^*E =
0$. We must show that $E = 0$. Let $T \in \DM^{gm}(k,A)$. It suffices to
prove that $\Hom(T,E) = 0$. Now by Proposition \ref{prop:bockstein-triangle} we
have the triangle $E \xrightarrow{\pi} E \to R\alpha_\pi^* L\alpha_{\pi\#} E =
0$. Thus multiplication by $\pi$ is an isomorphism on $\Hom(T, E)$. Let $K =
Frac(A)$. Since $\pi$ was arbitrary it follows that $\Hom(T,E)$ is a $K$-vector
space. Since $K \otimes_A K \ne 0$ we concude that $\Hom(T,E) = 0$ provided that
$\Hom(T,E) \otimes_A K = 0$. Let $\alpha_0: A \to K$ be the (flat) localisation.
By proposition \ref{prop:alpha-flat} we know that
$\Hom(T,E) \otimes_A K = \Hom(L\alpha_{0\#} T, L\alpha_{0\#} E)$, so it suffices to
show that $L\alpha_{0\#}E = 0$. But $K$ is of characteristic zero, so by
proposition \ref{prop:f*-conservative} it is enough to show that $f^*
L\alpha_{0\#} E = 0$. Since $L\alpha_{0\#}$ and $f^*$ ``commute'', this follows
from the assumption that $f^*E = 0$.

Now we prove Pic-injectivity. Let $E \in \DM(k,A)$ be such that $f^*E \approx
\tunit_{k^s}$ and $L\alpha_{\pi\#}E \approx \tunit_{A/(\pi)}$. As a first step,
I claim that there exists a finite extension $k \subset l \subset k^s$ such that
$g^*E \approx \tunit_l$, where $g: Spec(l) \to Spec(k)$. As in the proof of
Proposition \ref{prop:f*-conservative} we find
that $\Hom(\tunit_{k^s}, f^*E) = \colim_{k \subset l
\subset k^s} \Hom(\tunit_l, (l/k)^* E)$, where the colimit is over finite
subextensions. Hence there exist $l$ and an element $t \in \Hom(\tunit_l, g^*E)$
such that $(k^s/l)^*(t)$ is an isomorphism. The commutative diagram
\begin{equation*}
\begin{CD}
\Hom(\tunit_l, g^* E) @>>> \Hom(\tunit_{k^s}, f^* E) \approx A \\
@VVV                         @VVV                              \\
\Hom(\tunit_{l,A/(\pi)}, L\alpha_{\pi\#} g^*E) @>{\approx}>>
                      \Hom(\tunit_{k^s,A/(\pi)}, L\alpha_{\pi\#} f^*E) \approx A/(\pi)
\end{CD}
\end{equation*}
shows that $L\alpha_{\pi\#}(t)$ is an isomorphism. Thus by the first part
(conservativity), $t$ is an isomorphism.

Now we consider $\Hom(\tunit_k, E)$. From the Bockstein triangles and the
assumption $L\alpha_{\pi\#} \approx \tunit_{A/(\pi)}$ we get the exact sequences
\begin{gather*}
   \Hom(\tunit_{A/(\pi)}, L\alpha_{\pi\#}E[-1]) = 0 \to \Hom(\tunit_k, E)
   \xrightarrow{\pi} \Hom(\tunit_k, E) \\
   \to \Hom(\tunit_{A/(\pi)}, L\alpha_{\pi\#} E) \approx A/(\pi)
   \to \Hom(\tunit_k, E[1])
\end{gather*}
It follows that $\Hom(\tunit_k, E)$ is a torsion-free $A$-module (hence abelian
group). Thus by transfer it follows that $\Hom(\tunit_k, E) \to \Hom(\tunit_l,
g^*E) \approx A$ is injective. Let us denote the image by $I \subset A$. This is
a free $A$-module (of rank zero or one).

Since $\Hom(\tunit_l, g^*(E)[1]) = 0$ it follows by transfer
that $\Hom(\tunit_k, E[1])$ is $[l:k]$-torsion. Choosing $\pi$ of sufficiently large
characteristic, we find that $A/(\pi) \to \Hom(\tunit_k, E[1])$ is the zero map.
Thus $I = \Hom(\tunit_k, E) \ne 0$, i.e. $I \approx A$. It follows that
$\Hom(\tunit_k, E) \to \Hom(\tunit_{A/(\pi)}, L\alpha_{\pi\#}E) \approx A/(\pi)$
is surjective for each $\pi$.

Consider the commutative diagram
\begin{equation*}
\begin{CD}
\Hom(\tunit_k, E) @>>> \Hom(\tunit_{l}, g^* E) \approx A \\
@V{(*)}VV                         @V{(**)}VV                              \\
\Hom(\tunit_{A/(\pi)}, L\alpha_{\pi\#} E) @>{\approx}>>
                      \Hom(\tunit_{l,A/(\pi)}, L\alpha_{\pi\#} g^*E) \approx A/(\pi)
\end{CD}
\end{equation*}
The map (**) is the natural surjection and (*) is surjective as we just proved.
It follows that $I + (\pi) = A$ for each $\pi$ and so $I = A$. Thus there exists
$t' \in \Hom(\tunit_k, E)$ with $g^*(t') = t$ an isomorphism. Considering the diagram again one
finds that $L\alpha_{\pi \#}(t')$ is also an isomorphism. Thus $t'$ is an
isomorphism (by the first part, again) and we are done.
\end{proof}

We need two more auxiliary results. For the first, let $f: Spec(l) \to Spec(k)$ be a
Galois extension with group $G$. If $M \in \DM^{gm}(k, A)$ then the $A$-module
$\Hom(\tunit, f^*M) \approx \Hom(M(Spec(l)), M)$ has a natural action by $G$
(coming from automorphisms of $Spec(l)$). We denote this action by $\kappa_M: G
\to Aut(\Hom(\tunit, f^*M))$.

\begin{proposition} \label{prop:f*-pic}
Let $f: Spec(l) \to Spec(k)$ be (finite) Galois and $[l:k]$ invertible in $A$.
Then the above construction yields an injective
homomorphism
\[ \kappa: Ker(f^*: Pic(\DM^{gm}(k, A)) \to Pic(\DM^{gm}(l, A))) \to \Hom(Gal(l/k), A^\times). \]
\end{proposition}
\begin{proof}
Suppose that $M \in Pic(\DM^{gm}(k,A))$, $f^* M \simeq \tunit$ and let us show that $M \simeq \tunit$ if and
only if the action is trivial. Necessity is clear, we show sufficiency.

Independent of the assumptions on $[l:k]$ and $M$ I claim we have the following: if $t:
\tunit_L \to f^*M$ is any morphism, then $f^*(tr(t)): \tunit_L \to f^*M$ is the sum
of the conjugates under the $G$-action. Indeed the action on $\Hom_L(\tunit_L,
f^*M) \approx \Hom_k(\tunit_L, M)$ comes from premultiplication by elements of
$\Hom_k(\tunit_L, \tunit_L)$, whereas transfer comes from premultiplication with
the adjunction morphism. Thus to prove the claim we may assume that $M=\tunit$ and
$t=\id$, in which case the result follows from \cite[Exercise 1.11]{lecture-notes-mot-cohom}.

Thus reinstating our assumptions, let $t: \tunit_L \to f^*M$ be an isomorphism
and assume that the $G$-action is trivial. Then $tr(t/[l:k]): \tunit \to M$ is
an isomorphism since $f^*(tr(t/[l:k])) = t$ is, by Proposition
\ref{prop:f*-conservative}.

Finally we have to prove that $\kappa$ is a homomorphism, i.e. that $\kappa_{M
\otimes N} = \kappa_M \kappa_N$. For this let us denote the adjunction
isomorphism $\Hom_{\DM(k, A)}(M(l), T) \to \Hom_{\DM(l,A)}(\tunit, f^*T)$
generically by $ad$. One checks that given $f \in \Hom(M(l), M), g \in
\Hom(M(l), N)$ then $ad(f) \otimes ad(g) = ad((f \otimes g) \circ \alpha)$,
where $\alpha: M(l) \to M(l) \otimes M(l)$ is the map corresponding to $l
\otimes l \to l, a \otimes b \to ab$. Next observe that $\alpha$ is
$G$-equivariant if $G$ acts diagonally on $M(l) \otimes M(l)$. The result
follows.
\end{proof}

For the statement of the next result, we need $\DM(l, A)$ even if $l$ is not perfect.
It is explained in the next section what we mean by that. Under our assumptions
on $A$, it is equivalent to $\DM(l^p, A)$, where $l^p$ is the perfect closure of
$l$.

\begin{lemma} \label{lemm:base-change-trick}
Let $k$ be a perfect field, $X/k$ a smooth variety, $A$ a ring in which the
exponential characteristic of $k$ is invertible, and $M \in
\DM(k,A)$.

If for all $n \in \ZZ$ and all $x \in X$ (not necessarily closed) we have that
$\Hom_{\DM(x, A)}(\tunit\{n\}, M_x) = 0$, then also for all $n \in \ZZ$ we have
$\Hom_{\DM(k, A)}(MX\{n\}, M) = 0$.
\end{lemma}
\begin{proof}
We will prove the result by induction on $\dim{X}$. Thus in order to prove it
for $X$ we may assume that $\Hom_{\DM(k,A)}(MX'\{n'\}, M) = 0$
for every smooth, locally closed $X' \subset X$
with $\dim{X'} < \dim{X}$, and every $n' \in \ZZ$
(because the residue fields of $X'$ form a subset of
those of $X$). If $\dim{X} = 0$ then $X$ is a disjoint union of spectra of
fields, and the result is clear.

To prove the general case,
we may assume that $X$ is connected. Let $n \in \ZZ$ and $\alpha \in
\Hom(MX\{n\}, M)$. It suffices to show that $\alpha = 0$.
By considering the generic point and using
continuity \cite[Example 2.6(2)]{integral-mixed-motives} we conclude that there
exists a non-empty open subvariety $U \subset X$ such that $\alpha|_U =
0$. Let $Z = X \setminus U$.

If $Z$ is empty there is nothing to do. Otherwise there exists a non-empty,
smooth, connected open subvariety $U_1 \subset Z$,
since $k$ is perfect.

Let $Z' = Z \setminus U_1$, $U' = U \cup U_1 = X \setminus Z'$.
Then $U'$ is smooth open in $X$ and we have $X \setminus U' = Z'$,
which is strictly smaller than $Z$. We shall prove that
$\alpha|_{U'} = 0$. By repeating this argument with $U$ replaced by $U'$
(i.e. Noetherian induction on $Z$) it
will follow that $\alpha = 0$.

Note that $U_1 = U' \setminus U$ is closed in $U'$,
say of codimension $c$. Thus we get the distinguished Gysin triangle
\[ MU\{n\} \to MU'\{n\} \to MU_1\{n-c\}. \]
Now $\Hom(MU_1\{n-c\}, M) = 0$ by
the induction on dimension.
Thus $\Hom(MU'\{n\}, M) \to \Hom(MU\{n\},M)$ is injective. But
$(\alpha|_{U'})|_U = \alpha|_U = 0$ by assumption, so $\alpha|_{U'} = 0$.
\end{proof}

\section{Weight Structures and the Geometric Fixed Points Functors}
\label{sec:weights}

In this section, we will use Bondarko's theory of weight structures to construct
``generalised geometric fixed points functors'' and prove that they have good
properties. We shall fix a coefficient ring
$\FF$ on which an integer $e$ is invertible, and only work with fields of
exponential characteristic $e$.

We shall have to deal with $\DM(k, \FF)$ for $k$ an imperfect field. There is
now a fairly complete theory of $\DM(X, \FF)$ for Noetherian schemes over a
field of exponential characteristic $e$ (assumed invertible in $\FF$)
\cite{integral-mixed-motives}. It satisfies the six functors formalism, in
particular \emph{continuity}. We recall that if $k$ is an imperfect field
with perfect closure $k^p$, then the pullback $\DM(k, \FF) \to \DM(k^p, \FF)$ is
an equivalence of categories \cite[Proposition 8.1
(d)]{integral-mixed-motives}. This means that essentially all properties
known over perfect fields hold over imperfect fields as well. We also mention
that all of the categories $\DM(X, \FF)$ afford DG-enhancements. (This is well
known if $k$ is a perfect field and hence holds for $k$ any field by the
previous remark, and this is all we need. But it is actually clear that the constructions in
\cite{integral-mixed-motives} all yield DG categories.)

We shall work extensively in this section with weight structures
\cite{bondarko2010weight}, which we now review rapidly. Recall that given a
category $\mathcal{C}$ and a full subcategory $\mathcal{D} \subset \mathcal{C}$,
we call $\mathcal{D}$ \emph{Karoubi-closed in $\mathcal{C}$} if $\mathcal{D}$ is
closed under retracts \cite[p. 11]{bondarko2010weight}. In other words whenever
$X \in \mathcal{C}$ and $\id_X$ factorises through an object of $\mathcal{D}$,
then $X \in \mathcal{D}$. For example, if $\mathcal{C}$ is Karoubi-closed itself,
then $\mathcal{D}$ is Karoubi-closed in $\mathcal{C}$ if and only if
$\mathcal{D}$ is a Karoubi-closed category, and strictly full in $\mathcal{C}$.

Similarly, given a category $\mathcal{C}$ and a full subcategory $\mathcal{D}
\subset \mathcal{C}$, by the \emph{Karoubi-closure of $\mathcal{D}$ in
$\mathcal{C}$}
we mean the full subcategory of $\mathcal{C}$ spanned by all the objects which
are retracts of objects of $\mathcal{D}$. For example, if $\mathcal{D} \subset
\mathcal{C}$ is strictly full and $\mathcal{D}$ is a Karoubi-closed category,
then $\mathcal{D}$ is Karoubi-closed in $\mathcal{C}$.

\begin{definition} Let $\mathcal{C}$ be a
triangulated category and $\mathcal{C}^{w\ge 0}, \mathcal{C}^{w\le 0} \subset
\mathcal{C}$ two classes of objects. We call this a weight structure if the
following hold:
\begin{enumerate}[(i)]
\item $\mathcal{C}^{w\ge 0}, \mathcal{C}^{w \le 0}$ are additive and
  Karoubi-closed in $\mathcal{C}$.
\item $\mathcal{C}^{w \ge 0} \subset \mathcal{C}^{w \ge 0}[1],$
      $\mathcal{C}^{w \le 0}[1] \subset \mathcal{C}^{w \le 0}$
\item For $X \in \mathcal{C}^{w \ge 0}, Y \in \mathcal{C}^{w \le 0}$ we have
      $\Hom(X, Y[1]) = 0$.
\item For each $X \in \mathcal{C}$ there is a distinguished triangle
      \[ B[-1] \to X \to A \]
      with $B \in \mathcal{C}^{w \ge 0}$ and $A \in \mathcal{C}^{w \le 0}$.
\end{enumerate}
\end{definition}

These axioms look quite similar to those of a $t$-structure, but in practice
weight structures behave rather differently. We call a decomposition as in (iv)
a \emph{weight decomposition}. It is usually far from being unique. We put
$\mathcal{C}^{w \ge n} = \mathcal{C}^{w \ge 0}[-n]$ and $\mathcal{C}^{w \le n} =
\mathcal{C}^{w \le 0}[-n]$. We also write $\mathcal{C}^{w > n} = \mathcal{C}^{w
\ge n+1}$ etc. The intersection $\mathcal{C}^{w=0} := \mathcal{C}^{w \ge 0} \cap
\mathcal{C}^{w \le 0}$ is called the \emph{heart} of the weight structure.

A weight structure is called \emph{non-degenerate} if $\cap_n \mathcal{C}^{w \ge
n} = 0 = \cap_n \mathcal{C}^{w \le n}$. It is called \emph{bounded} if $\cup_n \mathcal{C}^{w \ge
n} = \mathcal{C} = \cup_n \mathcal{C}^{w \le n}$

A functor $F: \mathcal{C} \to \mathcal{D}$ between categories with weight
structures is called \emph{w-exact} if $F(\mathcal{C}^{w \le 0}) \subset
\mathcal{D}^{w \le 0}$ and $F(\mathcal{C}^{w \ge 0}) \subset
\mathcal{D}^{w \ge 0}$. It is called \emph{w-conservative} if given $X \in
\mathcal{C}$ with $F(X) \in \mathcal{D}^{w \le 0}$ we have $X \in \mathcal{C}^{w
\le 0}$, and similarly for $w \ge 0$. Note that a $w$-conservative
functor on a non-degenerate weight structure is conservative.

In the following proposition we summarise properties of weight structures we
use.

\begin{proposition} \label{prop:weights}
\begin{enumerate}[(1)]
\item $\mathcal{C}^{w \le 0}$ and $\mathcal{C}^{w \ge 0}$ are extension-stable:
  if $A \to B \to C$ is a distinguished triangle and $A, C \in \mathcal{C}^{w \le
  0}$ (respectively $A, C \in \mathcal{C}^{w \ge 0}$) then $B \in \mathcal{C}^{w
  \le 0}$ (respectively $B \in \mathcal{C}^{w \ge 0}$).

  Moreover $X \in \mathcal{C}^{w \ge 0}$ if and only if $\Hom(X, Y) = 0$ for all
  $Y \in \mathcal{C}^{w < 0}$, and similarly $X \in \mathcal{C}^{w \le 0}$ if
  and only if $\Hom(Y, X) = 0$ for all $Y \in \mathcal{C}^{w > 0}$.

\item Bounded weight structures are non-degenerate.

\item If $\mathcal{C}$ admits a DG-enhancement and the weight structure is
  bounded, then there exists a $w$-exact, $w$-conservative triangulated functor
  \[ t: \mathcal{C} \to K^b(\mathcal{C}^{w=0}) \]
  called the \emph{weight complex}. Its restriction to $\mathcal{C}^{w=0}$ is
  the natural inclusion.

\item If the weight structure is bounded and $\mathcal{C}^{w=0}$ is
  Karoubi-closed then so is $\mathcal{C}$.

\item If $H \subset \mathcal{C}$ is a \emph{negative} subcategory of a
  triangulated category (i.e. for $X, Y \in H$ we have $\Hom(X, Y[n]) = 0$ for
  $n>0$)
  generating it as a thick subcategory, then there exists a unique weight
  structure on $\mathcal{C}$ with $H \subset \mathcal{C}^{w=0}$. Moreover
  $\mathcal{C}^{w \le 0}$ is the smallest extension-stable Karoubi-closed
  subcategory of $\mathcal{C}$ containing $\cup_{n \ge 0} H[n]$, and similarly
  for $\mathcal{C}^{w \ge 0}$. The weight structure is bounded and
  $\mathcal{C}^{w=0}$ is the Karoubi-closure of $H$ in $\mathcal{C}$.

\item If $\mathcal{D} \subset \mathcal{C}$ is a triangulated subcategory such
  that $\mathcal{D}^{w\le 0} := \mathcal{D} \cap \mathcal{C}^{w\le 0}$ and
  $\mathcal{D}^{w\ge 0} := \mathcal{D} \cap \mathcal{C}^{w\ge 0}$ define a weight
  structure on $\mathcal{D}$ (we say the weight structure \emph{restricts} to
  $\mathcal{D}$) then the Verdier quotient $\mathcal{C}/\mathcal{D}$ affords a
  weight structure with $(\mathcal{C}/\mathcal{D})^{w\le 0}$ the Karoubi-closure
  of the image of $\mathcal{C}^{w \le 0}$ in $\mathcal{C}/\mathcal{D}$, and
  similarly for $(\mathcal{C}/\mathcal{D})^{w\ge 0}$,
  $(\mathcal{C}/\mathcal{D})^{w=0}$.

  The natural ``quotient'' functor $Q: \mathcal{C} \to \mathcal{C}/\mathcal{D}$
  is $w$-exact. If $X, Y \in \mathcal{C}^{w=0}$ then
  \[ \Hom(QX, QY) = \Hom(X, Y)/\Sigma_{Z \in \mathcal{D}^{w=0}} \Hom(Z, Y)
     \circ \Hom(X, Z). \]

  The weight structure on $\mathcal{C}/\mathcal{D}$ is bounded if the one on
  $\mathcal{C}$ is.
\end{enumerate}
\end{proposition}
\begin{proof}
(1) \cite[Proposition 1.3.3 (1-3)]{bondarko2010weight}.
(2) \cite[Proposition 1.3.6 (3) and comment after proof]{bondarko2010weight}.
(3) \cite[Proposition 3.3.1 (I), (IV) and Section 6.3]{bondarko2010weight}.
(4) \cite[Lemma 5.2.1]{bondarko2010weight}.
(5) \cite[Theorem 4.3.2 (II) and its proof]{bondarko2010weight}, \cite[Remark 2.1.2]{bondarko2016constructing}.
(6) \cite[Proposition 8.1.1]{bondarko2010weight}. Weight exactness holds by
definition of the weight structure on $\mathcal{C}/\mathcal{D}$.
\end{proof}

We shall call a triangulated category with a fixed weight structure a
$w$-category.

\begin{lemma} \label{lemm:sub-decomp}
Let $\mathcal{C}$ be a $w$-category with heart $H$, and $H' \subset H$ an
additive subcategory. Let $\mathcal{C}'$ be the thick triangulated subcategory generated
by $H'$ inside $\mathcal{C}$.

Then the weight structure of $\mathcal{C}$ restricts to $\mathcal{C}'$. In
particular, if $X \in \mathcal{C}'$ then we may
choose a weight decomposition $A \to X \to X'$ (i.e. $A \in \mathcal{C}^{w \ge
0}$ and $X' \in \mathcal{C}^{w < 0}$) with $A, X' \in \mathcal{C}'$.
\end{lemma}
\begin{proof}
This is just Proposition \ref{prop:weights} (5) which says that $\mathcal{C}'$,
being negatively generated by $H'$, carries a natural unique weight structure.
By the description provided we find $\mathcal{C}'^{w \le 0} \subset
\mathcal{C}^{w \le 0}, \mathcal{C}'^{w \ge 0} \subset \mathcal{C}^{w \ge 0}$.
Hence a weight decomposition in $\mathcal{C}'$ is also a weight decomposition in
$\mathcal{C}$. The rest follows from the definitions. (It follows from the
orthogonality characterisation that $\mathcal{C}'^{w \le 0} = \mathcal{C}^{w \le
0} \cap \mathcal{C'}$, but we do not need this.)
\end{proof}

\begin{lemma} \label{lemm:w-exact-functor}
Let $F: \mathcal{C} \to \mathcal{D}$ be a triangulated functor of
$w$-categories, and assume that the weight structure on $\mathcal{C}$ is
bounded. Then $F$ is $w$-exact if and only if $F(\mathcal{C}^{w=0}) \subset
\mathcal{D}^{w=0}$.
\end{lemma}
\begin{proof}
Necessity is clear, we show sufficiency.
We find by induction that the thick subcategory of $\mathcal{C}$ generated by $\mathcal{C}^{w=0}$
contains $\mathcal{C}^{w \le n} \cap \mathcal{C}^{w \ge -n}$ for all $n$, and
hence all of $\mathcal{C}$ by boundedness. It follows that the weight structure
on $\mathcal{C}$ is the one described in Proposition \ref{prop:weights} (5),
i.e. $\mathcal{C}^{w \ge 0}, \mathcal{C}^{w \le 0}$ are obtained as extension
closures of $\bigcup_{n\ge 0} \mathcal{C}^{w=n}, \bigcup_{n\le 0}
\mathcal{C}^{w=n}$. The result follows since $\mathcal{D}^{w \ge 0},
\mathcal{D}^{w \le 0}$ are extension-stable.
\end{proof}

\begin{lemma} \label{lemm:w-exact-tensor}
Let $\mathcal{C}$ be a $w$-category which is also a tensor category.
Assume that $\tunit_\mathcal{C} \in \mathcal{C}^{w=0}$ and that
tensoring is weight-bi-exact, i.e. that $\mathcal{C}^{w \le 0} \otimes
\mathcal{C}^{w \le 0} \subset \mathcal{C}^{w \le 0}$ and similarly for
$\mathcal{C}^{w \ge 0}$.

Then the weight complex functor is tensor
whenever $\mathcal{C}$ affords a tensor DG-enhancement and Pic-injective whenever
additionally the weight structure is bounded.

If moreover $\mathcal{C}$ is rigid
then the dualisation $D: \mathcal{C}^{op} \to \mathcal{C}$ is $w$-exact (i.e.
$D(\mathcal{C}^{w \ge 0}) \subset \mathcal{C}^{w \le 0}$ and vice versa).
\end{lemma}
\begin{proof}
If $\mathcal{D}$ is a negative DG tensor category, then $H^0(\mathcal{D})$ is
tensor in a natural way and the weight complex
functor $t$ manifestly respects the tensor structure. If $\mathcal{C}$ is a
tensor DG category with the property that $H^n(\Hom(X, Y)) = 0$ for all $X, Y
\in \mathcal{D}$ and $n > 0$ then the good truncation $\tau_{\le 0} \mathcal{D}$
is tensor in a natural way, and the quasi-equivalence $\tau_{\le 0} \mathcal{D}
\to \mathcal{D}$ is a tensor equivalence.

Hence the weight complex functor is
tensor as soon as there is any tensor DG enhancement of $\mathcal{C}^{w=0}$.
Moreover by Proposition
\ref{prop:weights} (3) if the weight structure is bounded then $t$ is
$w$-conservative. Since it induces an isomorphism on hearts it is a fortiori
Pic-injective. This proves the first part.

For the second part, let $X \in \mathcal{C}$. The category $\mathcal{C}$
being rigid means that there exists an object $DX$ such that $\otimes DX$ is
both right and left adjoint to $\otimes X$.

If $X \in \mathcal{C}^{w \ge 0}$ and $Y \in \mathcal{C}^{w > 0}$ then $\Hom(Y,
DX) = \Hom(Y \otimes X, \tunit) = 0$ because $Y \otimes X \in \mathcal{C}^{w >
0}$ whereas $\tunit \in \mathcal{C}^{w=0}$. It follows that $DX \in
\mathcal{C}^{w \le 0}$ by Proposition \ref{prop:weights} (1). The case of $X \in
\mathcal{C}^{w \le 0}$ is similar.
\end{proof}

We point out that for any field $k$, the category $\DM^{gm}(k, \FF)$ carries a
canonical weight structure \cite{bondarko2011}. (Note that the perfectness
assumption in that article can be dispensed with by passing to the equivalent
category $\DM^{gm}(k^p, \FF)$.) It is bounded and $\DM^{gm}(k, \FF)$ is also a
rigid tensor category with the tensor structure satisfying the assumptions of
Lemma \ref{lemm:w-exact-tensor}. The base change functors $f^*: \DM^{gm}(l, \FF)
\to \DM^{gm}(l', \FF)$ for $f: Spec(l')
\to Spec(l)$ are $w$-exact by Lemma \ref{lemm:w-exact-functor}.
The heart of the weight structure is $Chow(k^p,
\FF)$ which contains $Chow(k, \FF)$ as a full subcategory by Lemma
\ref{lemm:chow-insep-basechange}.

In the remainder of this section we will be dealing with the following
situation. The coefficient ring $\FF$ is a finite field of characteristic $p$
(necessarily $p \ne e$, where $e$ is the exponential characteristic of the
ground field $k$). For every extension $l/k$ we are given a set $S_l \subset
SmProj(l)$ such that for all closed points $x \in X \in S_l$ we have $p |
deg(x)$. Recall the categories $\langle S_l \rangle^{\otimes, T}_{Chow(l, \FF)}$
of Section \ref{sec:chow-motives}. We will assume that they are stable by base
change, i.e. that for $X \in S_l$ and $l'/l$ another extension we have $M X_{l'}
\in \langle S_{l'} \rangle^{\otimes, T}_{Chow(l', \FF)}$.

\newcommand{\DSTM}{\mathbf{D}\langle S \rangle\mathbf{TM}}

We write $\DSTM(l, \FF)$ for the thick triangulated subcategory
of $\DM(l, \FF)$ generated by $\langle S_l \rangle^{\otimes, T}_{Chow(l, \FF)}
\subset Chow(l, \FF) \subset \DM(l, \FF)^{w=0}$. It is tensor. The categories
$\DSTM(l, \FF)$ are also stable by base change in the sense that
if $f: Spec(l') \to Spec(l)$ is a field extension then $f^*(\DSTM(l, \FF))
\subset \DSTM(l', \FF)$. By Proposition \ref{prop:weights} (5) the weight
structure on $\DM^{gm}(l, \FF)$ restricts to $\DSTM(l, \FF)$, and the heart is
$\langle S_l \rangle^{\otimes, T}_{Chow(l, \FF)}$.

We write $\langle S_l \rangle^\otimes_{Chow(l, \FF)}$ for the Karoubi-closed
tensor subcategory of $Chow(l, \FF)$ generated by Tate twists of motives of varieties in $S_l$
(i.e. this is $\langle
S_l \rangle^{\otimes, T}_{Chow(l, \FF)}$ ``without the Tate motives''). By
Proposition \ref{prop:tatefree-presentation} this subcategory consists of
Tate-free objects. Let $\langle S_l\rangle^{tri} \subset \DSTM(l, \FF)$ be the
thick triangulated subcategory generated by $\langle S_l \rangle^{\otimes}_{Chow(l,
\FF)}$. As before, the weight structure restricts to $\langle S_l
\rangle^{tri}$. We write $\varphi^l_0: \DSTM(l, \FF) \to \DSTM(l, \FF)/\langle S
\rangle^{tri}$ for the Verdier quotient.

\begin{proposition} \label{prop:constructing-varphi}
The category $\DSTM(l, \FF)/\langle S \rangle^{tri}$ carries natural weight and
tensor structures, and $\varphi^l_0$ is a $w$-exact tensor functor. The
composite
\[ Tate(\FF) \to \DSTM(l, \FF) \to \DSTM(l, \FF)/\langle S \rangle^{tri} \]
is a full embedding with essential image
$\left( \DSTM(l, \FF)/\langle S \rangle^{tri} \right)^{w=0}$.
\end{proposition}
\begin{proof}
The existence of the weight structure and weight exactness is Proposition
\ref{prop:weights} (6). This also says that $\left( \DSTM(l, \FF)/\langle S
\rangle^{tri} \right)^{w=0}$ is generated as a Karoubi-closed subcategory by
$\varphi^l_0\left(\DSTM(l, \FF)^{w=0}\right)$. If $M \in \DSTM(l, \FF)^{w=0} =
\langle S_l \rangle^{\otimes, T}_{Chow(l, \FF)}$ then we may write $M \approx M'
\oplus T$ with $T$ a Tate and $M' \in \langle S_l \rangle^{\otimes}_{Chow(l,
\FF)}$, by Proposition \ref{prop:tatefree-presentation}. Thus $\varphi^l_0(M)
\approx \varphi^l_0(T)$ and so $\varphi^l_0: Tate(\FF) \to \left( \DSTM(l,
\FF)/\langle S \rangle^{tri} \right)^{w=0}$ is essentially surjective up to
(relative) Karoubi-closing. We shall show it is fully faithful whence its essential
image is (absolutely) Karoubi-closed and so $\varphi^l_0: Tate(\FF) \to \left( \DSTM(l,
\FF)/\langle S \rangle^{tri} \right)^{w=0}$ will be an equivalence. But by the
description in Proposition \ref{prop:weights} (6) it suffices to prove that any
morphism between Tate objects factoring through $\langle S_l
\rangle^\otimes_{Chow(l, \FF)}$ is zero. This follows from Lemma
\ref{lemm:tatefree-factoring}.

For the existence of the tensor structure we need
$\langle S \rangle^{tri} \otimes \DSTM(l, \FF) \subset \langle S \rangle^{tri}$;
then $\varphi^l_0$ is automatically tensor. Considering generators, it suffices
to show that $\langle S_l \rangle^{\otimes}_{Chow(l, \FF)} \otimes \langle S_l
\rangle^{\otimes, T}_{Chow(l, \FF)} \subset \langle S_l \rangle^\otimes_{Chow(l,
\FF)}$. This follows from Proposition \ref{prop:tatefree-presentation}.
\end{proof}

Let $l/k$ be any extension. We write $\Phi^l: \DSTM(k, \FF) \to K^b(Tate(\FF))$
for the composite
\begin{gather*}
\Phi^l: \DSTM(k, \FF) \to \DSTM(l, \FF) \to
     \DSTM(l, \FF)/\langle S_l \rangle^{tri} \\
  \xrightarrow{t}
     K^b\left( \left( \DSTM(l, \FF)/\langle S \rangle^{tri} \right)^{w=0} \right)
     \approx K^b(Tate(\FF))
\end{gather*}
of base change, the Verdier quotient functor $\varphi^l_0$, and the weight complex $t$.
It is a $w$-exact triangulated tensor functor. We can now state the main theorem of this section.

\begin{theorem}\label{thm:constructing-Phi}
Let $k$ be a ground field of exponential characteristic $e$, $\FF$ a finite
field of characteristic $p \ne e$. Suppose given for each field extension $l/k$
a set $S_l \subset SmProj(l)$ and a function $ex=ex_l: S_l \to \NN$. Assume that
the following hold (for all fields $l/k$):
\begin{enumerate}[(1)]
\item For $x \in X \in S_l$ closed, $p | deg(x)$.
\item If $l'/l$ is a field extension and $X \in S_l$ has no rational point over
      $l'$, then $X_{l'}$ is isomorphic to an object of $S_{l'}$ and
      $ex(X_{l'}) \le ex(X)$.
\item If $l'/l$ is a field extension and $X \in S_l$ has a rational point over
      $l'$, then $MX_{l'}$ is a summand of a motive of the form
 \[ T \oplus M(X_1^{(1)} \otimes \dots \otimes X_{n_1}^{(1)})\{i_1\} \oplus \dots
          \oplus M(X_1^{(m)} \otimes \dots \otimes X_{n_m}^{(m)})\{i_m\}, \]
      with $T \in Tate(\FF)$, $X_i^{(j)} \in S_{l'}$ and $ex(X_i^{(j)}) < ex(X)$ for
      all $i, j$.
\end{enumerate}
Then the family $\{\Phi^l\}_l$, as $l$ runs through finitely generated extensions
of $k$, is $w$-conservative (so in particular conservative) and Pic-injective.
\end{theorem}
We note that (2) and (3) imply that $\langle S_l \rangle^{\otimes, T}_{Chow(l,
\FF)}$ are stable by base change, i.e. we are in the situation we have been
discussing. Also (1) implies that none of the $X \in S_l$ have rational points
over $l$. The somewhat obscure functions $ex_l$ are necessary to make an
induction step in the proof work. We will always use $ex = \dim$ in
applications.

Before proving the result we explain how to compute $\Phi^l$ in the case that
$k$ is perfect (but $l$ need not be).
\begin{proposition} \label{prop:computing-Phi}
Assume in addition that $k$ is perfect. Let $l/k$ be a field extension.

There
exists an essentially unique additive functor $\Phi^l_0: \langle S_l \rangle^{\otimes,
T}_{Chow(l, \FF)} \to Tate(l, \FF)$ such that $\Phi^l_0|_{Tate(l, \FF)} = \id$ and
$\Phi^l_0(M) = 0$ if $M$ is Tate-free. It is tensor and the following diagram
commutes (up to natural isomorphism; the lower horizontal arrow is base change
of Chow motives):
\begin{equation*}
\begin{CD}
\DSTM(k, \FF) @>{\Phi^l}>>  K^b(Tate(\FF)) \\
@V{t}VV                     @A{\Phi^l_0}AA \\
K^b\left(\langle S_k \rangle^{\otimes, T}_{Chow(k, \FF)}\right)
  @>>> K^b\left(\langle S_l \rangle^{\otimes, T}_{Chow(l, \FF)}\right)
\end{CD}
\end{equation*}
\end{proposition}
\begin{proof}
Certainly $\Phi^l_0$ is essentially unique, using e.g. Proposition
\ref{prop:tatefree-presentation}. The functor $t \circ \varphi^l_0$ satisfies
the required properties, so $\Phi^l_0$ exists. It is tensor by construction.

To establish the commutativity claim, consider the diagram
\begin{equation*}
\begin{CD}
\DSTM(k, \FF) @>t>> K^b\left(\langle S_k \rangle^{\otimes, T}_{Chow(k, \FF)}\right) \\
@VVV                  @VVV \\
\DSTM(l, \FF) @>t>> K^b\left(\langle S_l \rangle^{\otimes, T}_{Chow(k, \FF)}\right) \\
@V{\varphi^l_0}VV     @V{\Phi^l_0}VV \\
\DSTM(l, \FF)/\langle S_l \rangle^{tri} @>{t}>> K^b(Tate(\FF)).
\end{CD}
\end{equation*}
It suffices to prove that the two squares commute (up to natural isomorphism).
This is most readily seen using DG-enhancements: let $\mathcal{D}(r)$ be a
functorial negative DG-enhancement of $\langle S_r \rangle^{\otimes,T}_{Chow(r, \FF)}
\subset \DSTM(r, \FF)$, for fields $r/k$. (In other words $\mathcal{D}(r)$ is a
DG-category with the same objects as $\langle S_r \rangle^{\otimes,T}_{Chow(r,
\FF)}$ and mapping complexes concentrated in non-positive degrees, such that these
mapping complexes compute morphisms in $\DSTM(r, \FF)$. Finally we ask that when
varying $r$, the assignment $r \mapsto \mathcal{D}(r)$ is a (pseudo-)functor and
$H^*(\mathcal{D}(r)) \to \DSTM(r, \FF)$ is a (pseudo-)natural transformation.)
Then it suffices to establish strict
commutativity of the diagram
\begin{equation*}
\begin{CD}
\mathcal{D}(k) @>>> \mathcal{D}(k)_0 \\
@VVV                  @VVV \\
\mathcal{D}(l) @>>> \mathcal{D}(l)_0 \approx \langle S_l \rangle^{\otimes, T}_{Chow(l, \FF)} \\
@VVV     @V{\Phi^l_0}VV \\
\mathcal{D}(l)/\langle S_l \rangle^{tri} @>>>
   \left( \mathcal{D}(l)/\langle S_l \rangle{tri} \right)_0 \approx Tate(\FF),
\end{CD}
\end{equation*}
where $\mathcal{D}_0$ for a negative DG-category means zero-truncation.
(Indeed the previous diagram is obtained by passing to
$Ho(\operatorname{Pre-Tr}(\bullet))$.) The upper square commutes by
definition and the lower square commutes if and only if it commutes on degree zero
morphisms, which is true essentially by definition of $\Phi^l_0$.
\end{proof}

We establish Theorem \ref{thm:constructing-Phi} through a series of lemmas.
\begin{lemma}
Let $\mathcal{C}$ be a $w$-category, $X \in \mathcal{C}^{w\le 0}$. Suppose given
weight decompositions $A \to X \to X'$ and $B[1] \to X' \to X''$ (i.e. $A, B \in
\mathcal{C}^{w \ge 0}$, $X' \in \mathcal{C}^{w<0}$ and $X'' \in \mathcal{C}^{w <
-1}$).

Then $A, B \in \mathcal{C}^{w=0}$ and for $T \in \mathcal{C}^{w=0}$ there is a
natural exact sequence
\[ \Hom(T, B) \to \Hom(T, A) \to \Hom(T, X) \to 0. \]
\end{lemma}
\begin{proof}
We have $A, B \in \mathcal{C}^{w=0}$ by (the dual of) \cite[Proposition 1.3.3
(6)]{bondarko2010weight}. There is an exact sequence
\[ \Hom(T, X'[-1]) \to \Hom(T, A) \to \Hom(T, X) \to \Hom(T, X') = 0 \]
where the last term is zero because $T \in \mathcal{C}^{w\ge 0}, X' \in \mathcal{C}^{w<0}$.
In particular $\Hom(T, A) \to \Hom(T, X)$ is surjective. Applying the same
reasoning to $\Hom(T, X'[-1])$ we find that $\Hom(T, B) \to \Hom(T, X'[-1])$ is
surjective and hence
\[ \Hom(T, B) \to \Hom(T, A) \to \Hom(T, X) \to 0 \]
is exact. This concludes the proof.
\end{proof}

\begin{corollary}
Let $X \in \DSTM(k, \FF)^{w \le 0}$ have a weight decomposition $T \to X \to X'$
with $T \in Tate(\FF)$ (and $X' \in \DSTM(k, \FF)^{w < 0}$). Suppose that
$\varphi^k(X) \in \left( \DSTM(k, \FF)/\langle S_k \rangle^{tri} \right)^{w<0}$.

Then for $T' \in Tate(\FF)$ we have $\Hom(T', X) = 0$.
\end{corollary}
\begin{proof}
Let $B[1] \to X' \to X''$ be a further weight decomposition.
Naturality in the above lemma yields the following commutative diagram with exact
rows
\begin{equation*}
\begin{CD}
\Hom(T', B) @>{\gamma}>> \Hom(T', T) @>>> \Hom(T', X) @>>> 0 \\
@V{\alpha}VV     @V{\beta}VV      @VVV      \\
\Hom(\varphi^k(T'), \varphi^k(B)) @>{\delta}>> \Hom(\varphi^k(T'), \varphi^k(T)) @>>> \Hom(\varphi^k(T'), \varphi^k(X)) @>>> 0.
\end{CD}
\end{equation*}
Since $\varphi^k$ is weight exact we have $\Hom(\varphi^k(T'), \varphi^k(X)) =
0$ and so $\delta$ is surjective. The construction of $\varphi^k$ (in particular
Proposition \ref{prop:constructing-varphi}) implies that $\alpha$ is surjective and $\beta$ is an isomorphism. It follows
that $\gamma$ is surjective, whence $\Hom(T', X) = 0$. This concludes the proof.
\end{proof}

The main work in proving our theorem is the following lemma. We let $\varphi^l:
\DSTM(k, \FF) \to \DSTM(l, \FF)/\langle S_l \rangle^{tri}$ be the composite of
$\varphi^l_0$ and base change.
\begin{lemma} \label{lemm:follows-1}
Let $X \in \DSTM(k, \FF)^{w \le 0}$ and suppose that for all $l/k$ finitely
generated, $\varphi^l(X) \in \left(\DSTM(l, \FF)/\langle S_l \rangle^{tri}
\right)^{w < 0}$. Then $X \in \DSTM(k, \FF)^{w < 0}$.
\end{lemma}
\begin{proof}
We begin by pointing out that Lemma \ref{lemm:base-change-trick} also applies if
$k$ is not perfect. Indeed if $k^p/k$ is the perfect closure then $X_{k^p}$ is
homeomorphic to $X$, so has the same set of points, and the residue field
extensions of $X_{k^p} \to X$ are purely inseparable, so induce equivalences on
$\DM(?, \FF)$. Thus the Lemma holds over $k$ if and only if it holds over $k^p$.

Let $\mathcal{R}$ be the set of finite multi-subsets of $\NN$ (i.e. the set of
finite non-increasing sequences in $\NN$). It is well-ordered lexicographically
and so can be used for induction. We extend $ex$ to a function $ex_l: \DSTM(l,
\FF) \to \mathcal{R}$. First, for $X_1, \dots, X_n \in S_l$ put $ex(X_1, \dots,
X_n) = \{\{ex(X_1), \dots, ex(X_n)\}\}$. Next, if $Y \in \DSTM(l, \FF)$ then
there exist $X_1, \dots, X_n \in S_l$ such that $Y \in \langle Tate(\FF), X_1,
\dots, X_n \rangle^{tri}$, i.e. $Y$ is in the thick tensor triangulated
subcategory generated by the $M X_i$ and the Tate motives.
We let $ex(Y)$ be the minimum of $ex(X_1, \dots, X_n)$ such that this holds. We
shall abuse notation and write $ex(Y) = ex(X_1, \dots, X_n)$ to additionally
mean that $Y \in \langle Tate(\FF), X_1, \dots, X_n\rangle^{tri}$.

Let us observe that if $ex(Y) = ex(X_1, \dots, X_n)$ and $l'/l$ is an extension
in which one of the $X_i$ acquires a rational point, then $ex(Y_{l'}) <
ex(Y_l)$, using assumptions (2) and (3).

We shall prove the result by induction on $ex(X)$. Note that is suffices to
prove that there is a weight decomposition $A \xrightarrow{\alpha} X \to X'$
(i.e. $A \in \DSTM(k, \FF)^{w=0}$ and $X' \in \DSTM(k, \FF)^{w<0}$) with
$\alpha=0$ (because then $X' \approx X \oplus A[1]$ and so $X \in
\DSTM(k,\FF)^{w<0}$, the latter being Karoubi-closed by definition).

If $ex(X) = \emptyset$ then $X$ must must be Tate. By Lemma
\ref{lemm:sub-decomp} we may choose a weight decomposition $T
\xrightarrow{\alpha} X \to X'$ with
$T \in Tate(\FF)$. By the corollary above (applied to $T'=T$) we find that
$\alpha=0$. This finishes the base case of our induction.

Suppose now $ex(X) = ex(X_1, \dots, X_n) > \emptyset$. If $l/k$ is any extension
such that one of the $X_1, \dots, X_n$ acquires a rational point over $l$, then
we may assume the lemma proved over $l$ by induction, so $X_l \in \DSTM(l,
\FF)^{w<0}$. Let $A \xrightarrow{\alpha} X \to X'$ be a weight decomposition; as before way may
choose $A \in \langle \{ X_1, \dots, X_n \} \rangle^{\otimes, T}_{Chow(k,
\FF)}$. Write $A \approx T \oplus A'$ as in Proposition
\ref{prop:tatefree-presentation}. I claim that $\alpha|_{A'} = 0$. It is enough
to show that if $Y$ is a product of the $X_i$ then $\Hom(MY\{n\},
X) = 0$ for all $n$. By Lemma \ref{lemm:base-change-trick}, it is enough to show
that for all $n \in \ZZ$ and $p \in Y$ we have that
$\Hom_{\DM(p, \FF)}(\tunit\{n\}, X_p) = 0$. But every variety has a rational
point after base change to any one of its points, so $X_p \in \DSTM(p,
\FF)^{w < 0}$ by induction. This proves the claim.

We thus have a weight decomposition $T \oplus A' \xrightarrow{(\alpha, 0)^T} X
\to X'$. Let $Y$ be a cone on $\alpha: T \to X$. We find that $X' \approx Y
\oplus A'[1]$ and hence $Y \in \DSTM(k, \FF)^{w<0}$. Thus $T
\xrightarrow{\alpha} X \to Y$ is a weight decomposition. Using the corollary
again we get $\Hom(T, X) = 0$ and so $\alpha = 0$. This finishes the induction
step.
\end{proof}

The rest of Theorem \ref{thm:constructing-Phi} is relatively easy to establish
now. We begin with the following.
\begin{lemma} \label{lemm:follows-2}
Let $\mathcal{C}, \mathcal{D}$ be $w$-categories with bi-$w$-exact tensor
structures. Suppose that $\mathcal{C}$ is rigid and its weight structure is
bounded.

Let $\phi: \mathcal{C} \to \mathcal{D}$ be a $w$-exact tensor functor such that
whenever $X \in \mathcal{C}^{w\le 0}$ and $\phi(X) \in \mathcal{D}^{w < 0}$ then
$X \in \mathcal{C}^{w<0}$.

Then $\phi$ is $w$-conservative.
\end{lemma}
\begin{proof}
Let $X \in \mathcal{C}$. If $\phi(X) \in \mathcal{D}^{w\le 0}$ then also $X \in
\mathcal{C}^{w\le 0}$. Indeed since the weight structure on $\mathcal{C}$ is
bounded we have $X \in \mathcal{C}^{w\le N}$ for some $N$. If $N > 0$ then the
assumptions imply that $X \in \mathcal{C}^{w \le N-1}$, and so on.

Suppose now instead that $\phi(X) \in \mathcal{D}^{w \ge 0}$. We need to show
that $X \in \mathcal{C}^{w \ge 0}$. But $X \in \mathcal{C}^{w \ge 0}$ if and
only if $DX \in \mathcal{C}^{w \le 0}$ by Lemma \ref{lemm:w-exact-tensor} (use
that $X \approx D(DX)$), and $\phi$ commutes with taking duals (since
$\mathcal{C}$ is rigid). Thus $\phi(DX) = D\phi(X) \in \mathcal{D}^{w \le 0}$,
so $DX \in \mathcal{C}^{w \le 0}$ and we are done.
\end{proof}
It follows from Lemmas \ref{lemm:follows-1} and \ref{lemm:follows-2} that
$\{\varphi^l\}_l$ is a $w$-conservative family. But all our
weight structures are bounded so the weight complex functors are
$w$-conservative, and thus $\{\Phi^l\}_l$ is also a $w$-conservative family.

Finally for $Pic$-injectivity, let $X \in \DSTM(k, \FF)$ be invertible with
$\Phi^l(X) \approx \tunit$ for all $l$. Since $\tunit \in K^b(Tate(\FF))^{w=0}$,
$w$-conservativity implies that $X \in \DSTM(k,\FF)^{w=0} = \langle S_k
\rangle^{\otimes, T}_{Chow(k, \FF)}$. Write $X \approx T \oplus X'$, with $T$
Tate and $X'$ Tate-free. Then $\tunit \approx \Phi^k(X) = T$ and so $T \approx
\tunit$. It follows that $\Phi^l(X) = \tunit \oplus \Phi^l(X') \in Tate(\FF)$. For this to be
invertible we need $\Phi^l(X') = 0$. Since this is true for all $l$,
conservativity implies that $X'=0$. This finishes the proof of Theorem
\ref{thm:constructing-Phi}.

\section{Application 1: Invertibililty of Affine Quadrics}
\label{sec:app1}

We now begin to reap in the benefits of the work of the previous sections. First
we construct the conservative and Pic-injective collection of functors we shall
use in the remainder of this work. After that we study invertibility of affine
quadrics.

We will be dealing with quadratic forms. If $l$ is a field and $\phi$ is a
non-degenerate quadratic form over $l$, we write $Y_\phi = Proj(\phi = 0)$ for
the projective quadric. This does not really make sense if $\dim{\phi} = 1$ in
which case we put $Y_\phi = \emptyset$ by convention.
Given $a \in l^\times$ we put $Y_\phi^a = Proj(\phi =
aZ^2)$ and $X_\phi^a = Spec(\phi = a)$. All
of these varieties are smooth.

Fix a perfect field $k$ of exponential characteristic $e \ne 2$
and coefficient ring $A$ containing $1/e$. We denote by $\QM(k,A)$ the
Karoubi-closed tensor subcategory of $Chow(k, A)$ generated by the (motives of) smooth
projective quadrics over $k$, and the Tate motives.

By \cite[Property (14.5.6)]{lecture-notes-mot-cohom} the category $Chow(k,A)$
embeds into $\DM^{gm}(k,A)$. We write $\DQM^{gm}(k,A)$ for the thick triangulated
subcategory of $\DM^{gm}(k,A)$ generated by $\QM(k,A)$. This is a tensor
category.

We write $\QM(k) = \QM(k,\ZZ[1/e])$ and $\DQM^{gm}(k) = \DQM^{gm}(k, \ZZ[1/e])$. As
promised, these categories contain the motives of all (smooth) affine quadrics.

\begin{lemma} \label{lemm:objects-in-DQM}
If $\phi$ is a non-degenerate quadratic form over the perfect field $k$ of
characteristic not two, and $a \in k^\times$, then
the affine quadric $X_\phi^a$ satisfies $M(X_\phi^a) \in \DQM^{gm}(k,A)$.
\end{lemma}
\begin{proof}
We have $X_\phi^a = Y_\phi^a \setminus Y_\phi$ and $M(Y_\phi^a), M(Y_\phi),
\tunit\{1\} \in \DQM^{gm}(k, A)$, so the result follows from the Gysin triangle.
\end{proof}

We recall the following result.
\begin{lemma}[(Rost)] \label{lemm:rost}
Let $\phi$ be an isotropic non-degenerate quadratic form. Then there exists a
non-degenerate form $\psi$ such that
\[M(Y_\phi) \approx \tunit \oplus M(Y_\psi)\{1\} \oplus \tunit\{\dim Y_\phi\}.  \]

Moreover for $a \in k^\times$ the natural ``inclusion'' $M(Y_\phi) \to M(Y_\phi^a)$
is given by
\[ \begin{pmatrix} \id & 0 & 0 \\ 0 & 0 & 0 \\ 0 & s\{1\} & i\{1\} \end{pmatrix}:
  \tunit \oplus \tunit\{\dim Y_\phi\} \oplus M(Y_\psi)\{1\}
\to \tunit \oplus \tunit\{\dim Y_\phi + 1\} \oplus M(Y_\psi^a)\{1\}, \]
where $i: M(Y_\psi) \to M(Y_\psi^a)$ is the natural ``inclusion'' and $s:
\tunit\{\dim Y_\psi^a\} \to M(Y_\psi^a)$ is the fundamental class (dual of the
structure map).
\end{lemma}
\begin{proof}
This is a result about Chow motives.

It is basically \cite[Proposition 2]{rost-pfister}. Rost starts with $\phi =
\mathbb{H} \perp \psi$, but this is equivalent to $\phi$ having a rational
point.

For the explicit form of the ``inclusion'', note first that all matrix entries shown
as zero have to be so for dimensional reasons. The entries ``$\id$'' and
``$i\{1\}$'' follow from naturality of Rost's construction. For the final entry,
we can argue as follows. Note that $\ZZ = CH^0(Y_\psi^a) = \Hom(\tunit\{\dim
Y_\psi^a+1\}, MY_\psi^a\{1\}) \approx \Hom(\tunit\{\dim
Y_\psi^a+1\}, MY_\phi) = CH^1(Y_\phi)$. The induced map we are interested in
corresponds under this identification to the cycle class of the closed
subvariety $Y_\phi \subset Y_\phi^a$. So up to verifying a sign (which is
irrelevant for all our applications), it is enough to
show that this class is a generator, which one sees for example by considering
the embedding into ambient projective space.
\end{proof}

\begin{lemma} \label{lemm:apply-thm-spl}
For a field extension $l/k$ let $S_l$ be the set of anisotropic projective
smooth quadrics over $l$, and let $ex_l: S_l \to \NN$ be the dimension function
$ex(X) = \dim{X}$. Then Theorem \ref{thm:constructing-Phi} applies, with
$\FF = \FF_2$.
\end{lemma}
We note that $\DSTM(k, \FF_2) = \DQM^{gm}(k, \FF_2)$, in the notation of the
Theorem.
\begin{proof}
Points on an anisotropic quadric have degree divisible by two by Springer's
theorem \cite[Chapter 7, Theorem 2.3]{lam-quadratic-forms}, hence condition (1)
holds. Condition (2) is satisfied essentially by definition. Finally condition
(3) follows from Lemma \ref{lemm:rost}.
\end{proof}

It follows from Lemma \ref{lemm:rost} that motives of quadrics are geometrically
Tate. Let $f: Spec(k^s) \to Spec(k)$ be a separable closure. It follows that the
weight complex functor $t:\DQM^{gm}(k^s) \to K^b(Chow(k^s, \ZZ[1/e]))$ takes
values in $K^b(Tate(\ZZ[1/e]))$. We write $\Psi$ for the composite
\[ \Psi: \DQM^{gm}(k) \xrightarrow{f^*} \DQM^{gm}(k^s) \xrightarrow{t} K^b(Tate(\ZZ[1/e])). \]

Let $g: Spec(l) \to Spec(k)$ be any field extension and $\alpha: \ZZ[1/e] \to \FF_2$
be the natural surjection. Via Lemma \ref{lemm:apply-thm-spl} and Theorem
\ref{thm:constructing-Phi} we obtain functors $\Phi^l: \DQM^{gm}(k, \FF_2) \to
K^b(Tate(\FF_2))$. We abuse notation and denote the composite with
change of coefficients
$\DQM^{gm}(k) \xrightarrow{L\alpha_\#} \DQM^{gm}(k, \FF_2) \to K^b(Tate(\FF_2))$
also by $\Phi^l$.

\begin{theorem} \label{thm:final-functors}
The functors $\Psi, \Phi^l$ are tensor triangulated. Together (as $l$ ranges
over all finitely generated extensions of $k$) they are conservative and
Pic-injective.
\end{theorem}
\begin{proof}
The functors are composites of tensor triangulated functors, so are tensor
triangulated.

By Theorem \ref{thm:conservative-on-DM} the collection $f^*,
\{L\alpha_{p\#}\}_p$ (where $p$ ranges over all primes) is
conservative and Pic-injective. Since all weight complex functors are
conservative and Pic-injective by Lemma \ref{lemm:w-exact-tensor}, the
collection $t f^*, \{t L\alpha_{p\#}\}_p$ is conservative and Pic-injective. We
have $tf^* = \Psi$. By Theorem \ref{thm:constructing-Phi} we may replace $t
L\alpha_{2\#}$ in our collection by $\{\Phi^l\}_l$.

It remains to deal with $L \alpha_{p\#}$ at odd $p$. Let $M \in \DQM^{gm}(k,
\ZZ[1/e])$. By repeated application of Lemma \ref{lemm:rost} we can find an
extension $L/k$ (which we may assume Galois) of degree a power of 2, such that
$M_L$ is in the triangulated subcategory generated by the Tate motives. In
particular $t(L\alpha_{p\#}M_L) \approx L\alpha_{p\#} \Psi(M)$ (as complexes of
Tate motives). Since $[L:k]$ is a power of two, base change along $L/k$ is
conservative in odd characteristic by Proposition \ref{prop:f*-conservative}.
Thus if $\Psi(M) \simeq 0$ then also $L\alpha_{p\#}M \simeq 0$ and our collection is
conservative.

We need to work a bit harder for Pic-injectivity. Let $M \in
\DQM^{gm}(k,\ZZ[1/e])$ be invertible with $\Phi^l(M) \simeq \tunit[0]$ for all
$l/k$ and $\Psi(M) \simeq \tunit[0]$. Then we know that $L\alpha_{2\#}(M) \simeq
\tunit$ by Theorem \ref{thm:constructing-Phi}. We also have $t(M_L) = \Psi(M)
\simeq \tunit$, so $M_L \simeq \tunit$ by Lemma \ref{lemm:w-exact-tensor}.
Consider the mod 2 Bockstein sequence
\begin{gather*}
  \Hom(\tunit, L\alpha_{2\#}M[-1]) = 0 \to \Hom(\tunit, M) \xrightarrow{2}
     \Hom(\tunit, M) \to \\
     \Hom(\tunit, L\alpha_{2\#}M) \to \Hom(\tunit, M[1])
     \xrightarrow{2} \Hom(\tunit, M[1]) \to \Hom(\tunit, L\alpha_{2\#}M[1]) = 0.
\end{gather*}
The extremal terms are zero because $L\alpha_{2\#}M \simeq \tunit$, and for the
same reason we have that $\Hom(\tunit, L\alpha_{2\#}M) \approx \FF_2$. Thus
$\Hom(\tunit, M)$ has no $2$-torsion, whereas $\Hom(\tunit, M[1])$ has no
$2$-cotorsion.
The composite $M \to M_L \to M$ of base change and transfer is multiplication by
$[L:k] = 2^N$. We conclude that $\Hom(\tunit, M)$ injects into $\Hom_L(\tunit_L,
M_L) \approx \ZZ[1/e]$ and that the kernel of $\Hom(\tunit, M[1]) \to
\Hom_L(\tunit_L, M_L[1]) = 0$ (i.e. the whole group) is contained in the
$2^N$-torsion. But multiplication by $2$ is surjective on $\Hom(\tunit, M[1])$,
whence so is multiplication by $2^N$, and we conclude that $\Hom(\tunit, M[1])
=0$.
Consequently we have $\Hom(\tunit, M) \approx \ZZ[1/e]$ (since it is an
ideal of $\ZZ[1/e]$ with a non-vanishing quotient, i.e. $\FF_2$).

We shall now apply Proposition \ref{prop:f*-pic}. As we have seen $M_L \simeq
\tunit$, so we obtain a $G=Gal(L/k)$-action on $\Hom(\tunit, M_L) \approx
\ZZ[1/e]$, i.e. a group homomorphism $\kappa_M: G \to \ZZ[1/e]^\times$. Since $e$ is prime
we have $\ZZ[1/e]^\times = \{\pm 1\} \times \{e^k | k \in \ZZ\}$ and since $G$ is
finite the image of $\kappa_M$ must be contained in $\{\pm 1\}$. Note that if
$\kappa_M = 1$ then $M \simeq \tunit$ and we are done. Indeed it suffices by Theorem
\ref{thm:conservative-on-DM} to show that $L \alpha_{p\#}M \simeq \tunit$ for
odd $p$. Since $(L\alpha_{p\#}M)_L \simeq \tunit$, by Proposition
\ref{prop:f*-pic} this happens if and only if an appropriate Galois action is
trivial, but this action is just the reduction
$G \xrightarrow{\kappa_M} \ZZ[1/e]^\times \to (\ZZ/p)^\times$.
So assume now that $\kappa_M$ is non-trivial.

Let $\beta:
\ZZ[1/e] \to \ZZ[1/(2e)]$ be the natural map. Note that $\kappa_M: G \to \{\pm 1\}$
has a kernel index 2, i.e. corresponds to a quadratic subextension $k
\subset k_2 \subset L$. I claim that $L\beta_\# M \simeq L\beta_\#
\tilde{M}Spec(k_2)$. Indeed this follows from Proposition \ref{prop:f*-pic}
applied to $A=\ZZ[1/(2e)]$, where $f^*$ becomes conservative, and
the observation that $\kappa_{\tilde{M}Spec(k_2)} = \kappa_M$.

In particular we must have $\Hom(\tunit, L\beta_\# \tilde{M}Spec(k_2)) \approx \Hom(\tunit, M)
\otimes_{\ZZ[1/e]} \ZZ[1/(2e)] = \ZZ[1/(2e)]$, by Proposition
\ref{prop:alpha-flat} and our previous computation. But one may compute easily
that $\Hom(\tunit, L\beta_\# \tilde{M}Spec(k_2)) = 0$. This contradiction concludes the proof.
\end{proof}

Note that if $A$ is a PID, then $Pic(K^b(Tate(A))) = \mathbb{Z} \oplus
\mathbb{Z}$. Consequently we have the following corollary.

\begin{corollary}
The abelian group $Pic(\DQM^{gm}(k, \ZZ[1/e]))$ is torsion-free (where $k$ is a
perfect field of exponential characteristic $e \ne 2$).
\end{corollary}

As the proof of the theorem shows, this is completely false with $\ZZ[1/(2e)]$
coefficients (or in the étale topology), where we have $Pic = \ZZ \oplus \ZZ
\oplus F$ where $F$ is an $\FF_2$-vector space.

We can now prove that affine quadrics are invertible.

\begin{theorem} \label{thm:invertibility}
Let $k$ be a perfect field of characteristic not two,
$\phi$ a non-degenerate quadratic form over $k$ and
$a \in k^\times$. Then $\tilde{M}(X_\phi^a)$ is invertible in $\DM^{gm}(k,
\ZZ[1/e])$.
\end{theorem}
\begin{proof}
We have $\tilde{M}X_\phi^a \in \DQM^{gm}(k)$ by
Lemma \ref{lemm:objects-in-DQM} and so we can use Theorem
\ref{thm:final-functors}. Since the category $\DQM^{gm}(k)$ is generated by
rigid objects (Chow motives) it is rigid and so conservative tensor functors detect invertibility, by standard
arguments. We thus need to show that $\Psi(\tilde{M}X_\phi^a)$ is
invertible and that for each $l/k$, $\Phi^l(\tilde{M}X_\phi^a)$ is invertible.

Let $d+2 = \dim{\phi}$. Let us put
$V_\phi^a = D(M X_\phi^a)\{d+1\}$ and $\tilde{V}_\phi^a = D(\tilde{M}
X_\phi^a)\{d+1\}$. Then $\tilde{M}X_\phi^a$ is invertible if and only if
$\tilde{V}_\phi^a$ is. From the closed inclusion $i: Y_\phi \to Y_\phi^a$ with complement
$X_\phi^a$ we get the dual Gysin triangle
\[ M Y_\phi \xrightarrow{i} M Y_\phi^a \to V_\phi^a. \]
It follows that $t(V_\phi^a) = [M Y_\phi \xrightarrow{i} \dot{M} Y_\phi^a]$.
Here the dot is used to indicate the term of degree zero in the chain complex.
Dualising the defining triangle of $\tilde{M} X_\phi^a$ we obtain
\[ \tunit\{d+1\} \xrightarrow{s} V_\phi^a \to \tilde{V}_\phi^a, \]
where $s$ is the fundamental class (dual of the structure map). Hence we finally
obtain
\[ t(\tilde V_\phi^a) = [M Y_\phi \oplus \tunit\{d+1\} \xrightarrow{(i, s)} \dot{M} Y_\phi^a] =: C(\phi). \]

The functor $\Psi$ is computed by first applying geometric base change, so
$\phi$ becomes completely split. In particular it has to be isotropic. An
induction on dimension using Lemma \ref{lemm:invertible-1} below shows that we
may reduce to $\dim{\phi} = 1$ or $2$, i.e. $\{x^2 = 1\}$ or $\{xy = 1\}$ (recall that
completely split quadrics are characterised by their dimension, so we can choose
any non-degenerate model quadric of the correct dimension). But
$\tilde{M}(\{x^2 = 1\}) = \tunit$ and $\tilde{M}(\{xy=1\})=\tilde{M}(\Gm)$ are both
invertible.

Dealing with $\Phi^l$ is a bit harder. 

The expression $C(\phi) \in K^b(\QM(k))$ makes sense even if $k$ is not perfect.
Using Proposition \ref{prop:computing-Phi} it suffices to prove:
if $l/k$ is any field extension, then $\Phi^l_0C(\phi_l)$ is invertible. We drop
the subscript zero from now on.
We may as well prove: if $k$ is any field and $\phi$ is any
non-degenerate quadratic form over $k$, then $\Phi^k(C(\phi))$ is invertible. By Lemma
\ref{lemm:invertible-1} below, if $\phi \approx \psi \perp \mathbb{H}$ then $C(\phi)
\simeq C(\psi)\{1\}$. We may thus assume that either $\phi$ is anisotropic, or $\phi
= \mathbb{H}$, or $\phi$ is of dimension one.

If $\phi = \mathbb{H}$ then $Y_\phi \approx Spec(k \times k)$, $Y_\phi^a
\approx \mathbb{P}^1$ and the result follows easily. If $\phi$ is of dimension
one then $MY_\phi = 0$ and either $MY_\phi^a = \tunit
\oplus \tunit$ or $M Y_\phi^a = M(k')$, where $k'/k$ is a quadratic extension.
Again the result follows easily.

So we may assume that $\phi$ is anisotropic. There are three cases. If $\phi \perp
\langle -a \rangle$ is also anisotropic, then none of $MY_\phi, MY_\phi^a$ afford
Tate summands, by Proposition \ref{prop:tatefree-presentation}. Thus
$\Phi^k(C(\phi)) = \tunit\{d+1\}[1]$ is invertible.

If
$\phi \perp \langle -a \rangle$ is isotropic, then $\phi \perp \langle -a
\rangle = \psi \perp \mathbb{H}$.
Suppose that $\psi$ has dimension greater than one. Then by (the contrapositive of)
Lemma \ref{lemm:invertible-2} below, $\psi$
is anisotropic. It follows that $MY_\phi^a \approx \tunit \oplus \tunit\{d+1\}
\oplus MY_\psi^a$ and $\Phi^k(C(\phi)) = [\tunit\{d+1\} \to \dot\tunit \oplus
\tunit\{d+1\}]$. The component $\tunit\{d+1\} \to \tunit\{d+1\}]$ comes from the
fundamental class of $M_\psi^a$ and so is an isomorphism. Thus $\Phi^k(C(\phi))
\simeq \tunit$ is invertible.

Finally it might be that $\psi$ has dimension one.
Then $Y_\phi^a \approx \mathbb{P}^1$ whereas $MY_\phi$ affords no Tate summands,
and the result follows as in the case of dimension greater than one.
This concludes the proof.
\end{proof}

\begin{lemma} \label{lemm:invertible-1}
Notation as in the theorem. If $\phi = \psi \perp \mathbb{H}$ then $C(\phi)
\simeq C(\psi)\{1\}$.
\end{lemma}
\begin{proof}
Using the explicit form for the inclusion $MY_\phi \to MY_\phi^a$ from Lemma
\ref{lemm:rost} we find that \[C(\phi) = [ (\tunit \oplus \tunit\{d\} \oplus
MY_\psi\{1\}) \oplus \tunit\{d+1\} \xrightarrow{\alpha} \dot\tunit \oplus \tunit\{d+1\} \oplus
MY_\psi^a\{1\} ],\]
where $\alpha$ is given by the matrix
\[ \begin{pmatrix}
 \id & 0 & 0 & 0 \\
  0  & 0 & 0 & f \\
  0  & s\{1\} & i\{1\} & 0
\end{pmatrix}. \]
Here $f$ comes from the fundamental class and so is an isomorphism. It follows
that $C(\phi) \approx C(\psi)\{1\} \oplus cone(\id_\tunit)[-1] \oplus
cone(\id_{\tunit\{d+1\}})[-1] \simeq C(\psi)\{1\}$. This is the desired result.
\end{proof}

\begin{lemma} \label{lemm:invertible-2}
If $\phi \perp \langle a \rangle \approx \psi \perp \mathbb{H} \perp
\mathbb{H}$, then $\phi$ is isotropic.
\end{lemma}
\begin{proof}
Let $X = Y_{\langle a \rangle \perp \phi}$. Then $Y_\phi = X \cap \{X_0 = 0\}$.
Since $\langle a \rangle \perp \phi \approx \psi \perp \mathbb{H} \perp \mathbb{H}$,
we find that $Y_{\mathbb{H} \perp \mathbb{H}} \subset X$. Then $Y_\phi \cap
Y_{\mathbb{H} \perp \mathbb{H}} = Y_{\mathbb{H} \perp \mathbb{H}} \cap \{X_0 =
0\}$ (intersecting inside $X$). Now we know that after a \emph{linear} change of
coordinates $(X_0: \dots: X_r) \mapsto (T_0: \dots : T_r)$ the subvariety
$Y_{\mathbb{H} \perp \mathbb{H}}$ of $X$ is given by the equations $T_0 T_1 + T_2
T_3 = 0$, $T_i=0$ for $i > 3$. Thus $Y_\phi \cap Y_{\mathbb{H} \perp
\mathbb{H}}$ is obtained by adding a further \emph{linear}
constraint in the $T_0, T_1, T_2, T_3$. It is easy to see that there must be a
rational, non-zero solution, so $Y_\phi$ has a rational point. This was to be shown.
\end{proof}

\section{Application 2: Po Hu's Conjectures for Motives}
\label{sec:app2}

In this final section we prove a version for motives of Po Hu's conjectures
\cite[Conjecture 1.4]{HuPicard}. We retain notation from the previous section.

For $\ul a = (a_1, \dots, a_n) \in (k^\times)^n$, $b \in k^\times$ let us put
\[ U_{\ul a}^b = X_{\langle \langle a_1, \dots, a_n \rangle \rangle}^b, \]
where $\langle \langle a_1, \dots, a_n \rangle \rangle$ is the $n$-fold Pfister
quadric associated with the symbol $\ul a$. We use notation such as $\ul{a}, a'
= (a_1, \dots, a_n, a') \in (k^\times)^{n+1}$ for concatenation of tuples.

\begin{theorem} \label{thm:hu-conj}
Let $k$ be a perfect field of characteristic not two, and $\ul a \in
(k^\times)^n, b \in k^\times$.

In $\DM^{gm}(k)$ there is an isomorphism
\begin{equation} \label{eq:hu-iso}
\tilde{M}(U_{\ul{a}, b}^1) \otimes \tilde{M}(U_{\ul a}^b)[1] \approx
\tilde{M}(U_{\ul a}^1)\{2^n\}.
\end{equation}
\end{theorem}

To prove this, we have to recall some facts about \emph{Rost motives}. If
$\ul a \in (k^\times)^n$, then there is the associated Rost motive $R_{\ul a}
\in \QM(k)$. Recall that one has $H^1_{et}(k, \FF_2) = k^\times/2,$ and hence
cup product yields a natural map $\partial = \partial^k: (k^\times)^n \to
H^n_{et}(k, \ZZ/2)$. The Rost motives have the remarkable property that $R_{\ul
a}$ is irreducible if and only if $\partial(\ul a) \ne 0$. In fact there are
canonical maps
\begin{equation} \label{eq:rost-triangle}
  \tunit\{2^{n-1} - 1\} \to R_{\ul a} \to \tunit
\end{equation}
(which we call structure maps) and if $\partial(\ul a) = 0$ then this is a
splitting distinguished triangle. The same statements hold true with $\FF_2$
coefficients. These results follow from the work of a number of people, see
\cite{morel-voevodskys-proof} for an overview.

The relationship between Rost motives and $U_{\ul a}^b$ is encapsulated in the
following proposition.

\begin{proposition} \label{prop:rost-U}
For $\ul a \in (k^\times)^n, b \in k^\times$ there is a distinguished triangle
\[ \tilde{M}(U_{\ul a}^b) \to R_{\ul a, b} \to R_{\ul a}\{2^{n-1}\} \oplus \tunit. \]
Here $R_{\ul a, b} \to \tunit$ is the structure map, and the composite
\[ \tunit\{2^n - 1\} \to R_{\ul a, b} \to R_{\ul a}\{2^{n-1}\} \]
is the $\{2^{n-1}\}$ twist of the structure map $\tunit\{2^{n-1}-1\} \to R_{\ul
a}$.
\end{proposition}
\begin{proof}
This is essentially \cite[proof of Proposition 5.5]{HuRemarks}.

We know that $U := U_{\ul a}^b$ is the complement of $X := Y_{\langle \langle
\ul a \rangle\rangle}$ in $Y := Y_{\langle \langle \ul a \rangle \rangle}^b.$ By the work
of Rost \cite[Theorem 17 and Proposition 19]{rost-pfister},
if we put $R_n := R_{\underline a, b}$ and $R_{n-1} = R_{\underline a},$ then
\[ M(Y) = R_n \oplus \bigoplus_{k=1}^{2^{n-1}-1} R_{n-1}\{k\} := R_n \oplus R',
  \quad M(X) = \bigoplus_{k=0}^{2^{n-1}-1} R_{n-1}\{k\} := R_{n-1} \oplus R' \]
and the natural map $M(X) \to M(Y)$ is the identity on $R'.$

The localisation triangle $M^c(X) = M(X) \to M^c(Y) = M(Y) \to M^c(U)$ fits into
the following commutative diagram of (distinguished) triangles:
\begin{equation*}
\begin{CD}
 R'   @=  R'             \\
@VVV       @VVV      @.  \\
M(X) @>>> M(Y) @>>> M^c(U) \\
@VVV      @VVV      @.   \\
R_{n-1} @. R_n           \\
\end{CD}
\end{equation*}
An application of the octahedral axiom yields a distinguished triangle $R_{n-1}
\to R_n \to M^c(U).$ Noting that $DM^c(U) = M(U)\{-(2^n-1)\}, DR_n =
R_n\{-(2^n-1)\}$ and $DR_{n-1} = R_{n-1}\{-(2^{n-1}-1\},$ by dualising and
twisting the triangle, we find a distinguished triangle $M(U) \to R_n \to
R_{n-1}\{2^{n-1}\}.$ Adding in the copy of $\tunit$ implied in $\tilde{M}(U),$ we
get the claimed triangle with the correct map $R_n \to \tunit$.

To see the second claim about the differential, the important point is that in the triangle $R_{n-1} \to R_n
\to M^c(U)$ the map $R_{n-1} \to R_n$ is induced from the inclusion $M(X) \to
M(Y)$ by passing to the appropriate summands. It follows that $R_{n-1}
\to R_n \to \tunit$ is the structure map of $R_{n-1} \to \tunit.$ The desired result
now follows by dualising.
\end{proof}

\begin{proof}[Proof of Theorem \ref{thm:hu-conj}]
By Lemma \ref{lemm:objects-in-DQM}, we have $\tilde{M}(U_{\ul a}^b) \in
\DQM^{gm}(k)$, etc. We also know by Theorem \ref{thm:invertibility} that both
sides of equation \eqref{eq:hu-iso} are invertible. Hence if $F: \DQM^{gm}(k)
\to \mathcal{C}$ is a Pic-injective functor, it suffices to prove that $F(LHS)
\approx F(RHS)$.

Of course we use the Pic-injective collection from Theorem
\ref{thm:final-functors}.

From Proposition \ref{prop:rost-U} we know that
\[ t(\tilde{M}(U_{\ul a}^b) = [ \dot{R}_{\ul a, b} \to R_{\ul a}\{2^n\} \oplus \tunit], \]
and we also know certain things about the differential. To compute $\Psi$, we
have to consider geometric base change, where the triangle
\eqref{eq:rost-triangle} is splitting distinguished. One obtains
\[ \Psi(\tilde{M}(U_{\ul a}^b)) = [\dot\tunit \oplus \tunit\{2^n - 1\} \to
         \tunit\{2^{n-1}\} \oplus \tunit\{2^n - 1\} \oplus \tunit] \]
and from the information about the differential given in proposition
\ref{prop:rost-U} we deduce that $\Psi(\tilde{M}(U_{\ul a}^b) \simeq
\tunit\{2^{n-1}\}[-1]$. Thus $\Psi(LHS) \approx \Psi(RHS)$ reads
\[ \tunit\{2^n\}[-1] \otimes \tunit\{2^{n-1}\}[-1][1] \approx
\tunit\{2^{n-1}\}[-1]\{2^n\}, \]
which is certainly true.

Now let $l/k$ be an arbitrary field extension. We need to prove $\Phi^l(LHS)
\approx \Phi^l(RHS)$. This involves $R_{\ul a}, R_{\ul a,b}, R_{\ul a, 1}$ and
$R_{\ul a, b, 1}$. Depending on $l$ these may or may not split into Tate
motives, so may or may not survive $\Phi$. We see that $R_{\ul a, 1}$ and
$R_{\ul a, b, 1}$ always split (because $\partial^l(1) = 0$), and that $R_{\ul
a,b}$ splits whenever $R_{\ul a}$ splits (because $\partial(\ul a, b) =
\partial(\ul a) \cup \partial(b)$).

If $R_{\ul a}$ splits then everything is split and $\Phi^l$ is just mod two
reduction of $\Psi$, so we know the equation is satisfied. Thus there are just
two cases and three things in each to compute, which we gather in Table \ref{tab:1}.

\begin{table}
\center
\caption{Terms needed to compute $\Phi^l$.}
\label{tab:1}
\begin{tabular}{c|c|c}
      & $\partial^l(\ul a, b) \ne 0$ & $\partial^l(\ul a, b) = 0$ but $\partial^l(\ul a) \ne 0$ \\
\hline
$\Phi^l(U_{\ul a, b}^1)$ & $[\dot\tunit \oplus \tunit\{2^{n+1}-1\} \to \tunit]$
   & $[\dot\tunit \oplus \tunit\{2^{n+1}-1\} \to \tunit\{2^n\} \oplus \tunit\{2^{n+1}-1\} \oplus \tunit]$ \\
$\Phi^l(U_{\ul a}^b)$    & $[\dot 0 \to \tunit]$ 
   & $[\dot\tunit \oplus \tunit\{2^n-1\} \to \tunit]$ \\
$\Phi^l(U_{\ul a}^1)$    & $[\dot\tunit \oplus \tunit\{2^n-1\} \to \tunit]$
   & $[\dot\tunit \oplus \tunit\{2^n-1\} \to \tunit]$ \\
\end{tabular}
\end{table}

The differentials can again
be figured out using Proposition \ref{prop:rost-U}. Using these one can simplify
the expressions. We have gathered the results in Table \ref{tab:2}.

\begin{table}
\center
\caption{Terms needed to compute $\Phi^l$, simplified form.}
\label{tab:2}
\begin{tabular}{c|c|c}
      & $\partial^l(\ul a, b) \ne 0$ & $\partial^l(\ul a, b) = 0$ but $\partial^l(\ul a) \ne 0$ \\
\hline
$\Phi^l(U_{\ul a, b}^1)$ & $\tunit\{2^{n+1} - 1\}$ & $\tunit\{2^n\}[-1]$ \\
$\Phi^l(U_{\ul a}^b)$    & $\tunit[-1]$   & $\tunit\{2^n-1\}$ \\
$\Phi^l(U_{\ul a}^1)$    & $\tunit\{2^n-1\}$  & $\tunit\{2^n-1\}$ \\
\end{tabular}
\end{table}

To complete the proof, we check that $\Phi^l(LHS) \approx \Phi^l(RHS)$ in both
cases. This is easy.
\end{proof}

\bibliographystyle{amsplain}
\bibliography{bibliography}

\Addresses
\end{document}